\documentclass{ws-m3as}
\usepackage[utf8x]{inputenc} 
\usepackage{amsmath,stmaryrd}
\usepackage{graphics}
\usepackage{color}

\title{A DRIFT-DIFFUSION-REACTION MODEL FOR EXCITONIC PHOTOVOLTAIC BILAYERS: ASYMPTOTIC ANALYSIS AND A 2-D HDG FINITE-ELEMENT SCHEME}
\author{DANIEL BRINKMAN}

\address{Department of Applied Mathematics and Theoretical Physics, University of Cambridge \\
 Centre for Mathematical Sciences, Wilberforce Road \\
Cambridge, CB3 0WA, United Kingdom \\
D.Brinkman@damtp.cam.ac.uk}

\author{KLEMENS FELLNER}

\address{Institute for Mathematics and Scientific Computing, University of Graz Heinrichgasse 36 \\
 8010 Graz, Austria\\
Klemens.Fellner@uni-graz.at}

\author{PETER A. MARKOWICH}

\address{Mathematical and Computer Sciences and Engineering Division \\
King Abdullah University of Science \& Technology \\
Thuwal 23955-6900, Kingdom of Saudi Arabia \\
P.A.Markowich@damtp.cam.ac.uk}

\author{MARIE-THERESE WOLFRAM}

\address{Department of Applied Mathematics and Theoretical Physics, University of Cambridge \\
Department of Mathematics, University of Vienna Nordbergstrasse 15 \\
1090 Vienna, Austria \\
marie-therese.wolfram@univie.ac.at}

\def\eps{\varepsilon}

\begin{document}
\maketitle
\begin{abstract}
We present and discuss a mathematical model for the operation of bilayer organic photovoltaic devices. Our model 
couples drift-diffusion-recombination equations for the charge carriers (specifically, electrons and holes) with a reaction-diffusion equation for the excitons/polaron pairs and Poisson's equation for the self-consistent electrostatic potential. The material difference (i.e. the HOMO/LUMO gap) of the two organic substrates forming the bilayer device are included as a work-function potential. 

Firstly, we perform an asymptotic analysis of the scaled one-dimensional stationary state system 
i) with focus on the dynamics on the interface and ii) with the goal of simplifying the bulk dynamics away for the interface. 
Secondly, we present a two-dimensional Hybrid Discontinuous Galerkin Finite Element numerical scheme which 
is very well suited to resolve i) the material changes ii) the resulting strong variation over the interface and iii) the necessary upwinding in the discretization of drift-diffusion equations. 
Finally, we compare the numerical results with the approximating asymptotics. 
\end{abstract}



\section{Introduction}
The search for cheap, environmentally sustainable energy solutions has lead to intense interest and research in the area of organic photovoltaic (OPV) devices. Currently, these devices have solar efficiencies which are significantly below those of modern inorganic semiconductor devices (i.e. at present 8.3\% compared to 27.6\% for ``one-sun'' Gallium Arsenide and 42.3\% for concentrator cells\cite{GEHW}), currently limiting commercial implementation. 

A main difference of OPV devices (in contrast to inorganic PV) is the generation mechanism of free charge carriers, which typically occurs through the generation and dissociation of so called excitons, excited energy states created by incoming light. (We shall discuss excitons and the generation of free charge carriers in more detail below.)

Many mathematical models have been proposed to study the behavior of bilayer OPV cells as simple implementable organic devices: Refs. \refcite{BC,Cetal,BRG,K}. These models involve various approximations of the generation of free charges. A mathematical basis for the these models comes from extensive literature on inorganic semiconductor models. For these devices, the standard macroscopic models couple two drift-diffusion-reaction equations for the behavior of the free charge carrier densities (i.e. electrons and holes) in a system with the Poisson equation governing the self-consistent electrical potential.\cite{M} Much is known about such systems, including existence and uniqueness results for the steady-state (see e.g. Refs. \refcite{M,Mss} and the references therein). The steady state is of particular interest for photovoltaic devices which are expected to generate power on time scales much longer than the duration of the transient (i.e. the switching-on) dynamics. 

Drift-diffusion-reaction type models for OPV devices must crucially take into account the specific material properties of inorganic and organic semiconductor materials, which involves in particular electric field dependent mobilities and specific recombination and dissociation rates (see Sec. \ref{sec:scaling} below). A second crucial amendment has to model the exciton dynamics and couple it appropriately to the original system.
 
We shall thus consider the following evolution equations for the four main components, i.e. the concentrations of electrons ($n$), holes ($p$), and excitons ($X$) and the electric potential in the device ($V$), which must be calculated self-consistently.

The evolution of the free charge carrier densities $n$ and $p$ shall be described by typical drift-diffusion-recombination equations and an additional source term due to the excitons:
\begin{align}
q\frac{\partial n}{\partial t} & = \nabla \cdot (q D_n \nabla n - q \mu_n n \nabla (U+V)) - q R_{np} + q k_d X \label{eq:basicelectrons} \\
q\frac{\partial p}{\partial t} &= \nabla \cdot (q D_p \nabla p + q \mu_p p \nabla (U+V)) - q R_{np} + q k_d X \label{eq:basicholes} \text{.}
\end{align}
Here $q$ denotes the positive fundamental charge, $D_n$, $D_p$ and $\mu_n$, $\mu_p$ are the diffusion coefficients and mobilities of $n$ and $p$, respectively, $R_{np}$ is the recombination rate of free electrons and free holes, and $k_d$ is the rate of dissociation of excitons into a pair consisting of a free electron and hole. Functional expressions for diffusion coefficients, mobilities and recombination/dissociation rates  in organic semiconductor materials will be specified in Sec. \ref{sec:scaling}.

\begin{figure}[htbp]
\centering
\scalebox{0.35}{\includegraphics{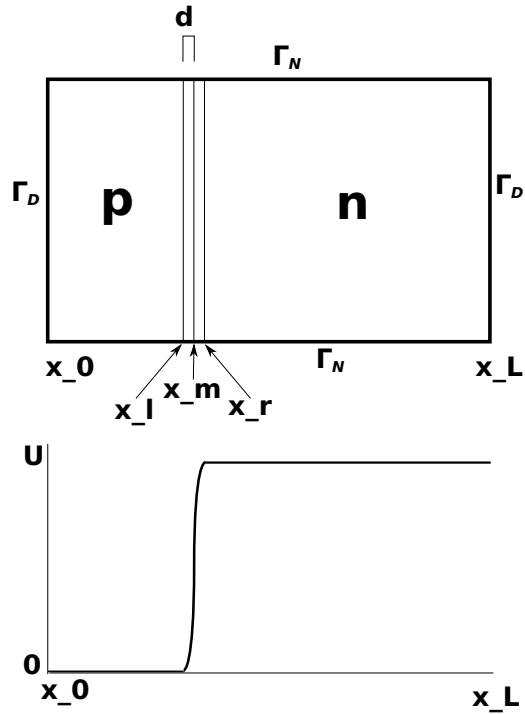}}
\caption{Schematic of the 2-D device giving labels to the quantities used in the asymptotics and a graphical representation of the function $U$. For the numerics we take $U$ to be a piecewise quadratic connecting the values at either side so that $U \in C^1$. Note that for the asymptotics the vertical direction of the first diagram is suppressed and the problem becomes 1-D.} \label{fig:schematic}
\end{figure}

For a bilayer OPV device formed of two different organic semiconductor materials (see Fig. \ref{fig:schematic}), the potential $U$ (which is not an electric potential and doesn't contribute to the electric field) models an additional convection term that comes from the differences in the energy levels of the two materials. In an organic device, the energy levels involved in charge transfer are the Highest Occupied Molecular Orbital (HOMO) and the Lowest Unoccupied Molecular Orbital (LUMO). The HOMO roughly corresponds to the valence band for classical semiconductors, and the LUMO resembles the conduction band. In general, the HOMO-LUMO gap corresponds to the band-gap. 

For most organic materials, the HOMO-LUMO gap is too large for a photon to create free electrons and holes. Instead, it creates a bound electron/hole pair, a so-called exciton. In the bulk material, these excitons usually recombine without producing free carriers (with rate $k_r$, see below). However, at an interface between two materials with suitable differences in the HOMO/LUMO properties, the excitons tend to align and split over the interface so that the electrons and holes are in separate (energetically favorable) materials. This effect is included in the model equations \eqref{eq:basicelectrons} and \eqref{eq:basicholes} by the term $U$, or more precisely, by the change in $U$  over the interface region.  

Excitons are typically called polaron pairs (also referred to as charge-transfer states or coulomb bound pairs) when aligned across the interface and are much more stable with respect to recombination. On the other hand, the aligned polaron pairs will either be pushed together or pulled apart by an approximately parallel electric field. Accordingly the dissociation rate ($k_d$) will be heavily field dependent. The combination of these two effects (reduced polaron pair recombination and field-driven dissociation) can lead to very high quantum efficiency (proportion of light creating free charges) for polaron pairs at an interface under an appropriately applied external potential.\cite{AHB} Note that this is an internal quantum efficiency for polaron pair dissociation, not an external quantum efficiency, which further includes how many of the incoming photons actually create polaron pairs.

Following the above discussion we shall assume that every exciton becomes immediately a polaron pair at the interface (see further explanation below). We therefore consider the following diffusion equation for the evolution of a combined density of excitons and polaron pairs (to be identified by their position), which lacks a convection term since the excitons are electrically neutral:
\begin{equation}
 \frac{\partial X}{\partial t} = \nabla \cdot (D_X \nabla X) + c R_{np} + G - k_d X - k_r X \text{.} \label{eq:basicexcitons}
\end{equation}
The coefficient $G$ is the photo-generation rate for excitons, $k_r$ is the rate of geminate recombination of the excitons, and $c$ represents the proportion of recombining free charges that form excitonic states (as opposed to recombination leading to emission of light, etc). 

Finally, we have the Poisson equation for the self-consistent electric potential:
\begin{equation}
  \epsilon_0 \nabla \cdot (\epsilon_r \nabla V) = q (n-p- h \chi_I \nabla_\nu X)  \text{,} \label{eq:basicpoisson} \\
\end{equation}
where $\epsilon_0$ denotes the vacuum permittivity and $\epsilon_r$ is the relative permittivity of the material. Furthermore, $\nabla_\nu$ represents the derivative normal to the direction of the material interface $I$ and $\chi_I$ denotes the indicator function of the interface region while $h$ is the separation length of the charge pairs of polaron pairs in the interface region. In fact,  the extra term proportional to $\nabla_\nu X$ is the field contribution due to the alignment of the polaron pairs at the material interface. The gradient results from taking the continuum limit of a sum of Dirac delta primes (corresponding to point dipoles). The exciton distribution $X$ includes both polaron pairs and excitons, with their identity based on their location in the system. Polaron-pairs are split over the interface such that the electron and hole (although still bound) each lie in an energetically favorable material. This has the dual effect of aligning a sheet of dipoles (resulting in the aforementioned term in Eq. \eqref{eq:basicpoisson}) and of making polaron pair states extremely energetically favorable to excitons.\cite{AHB} We thus assume that any $X$ in an interface region $I$ are polaron pairs and that any $X \in I^C$ (the complement of $I$) are excitons. In actual devices, organic interfaces are not sharp, and can be blended over a variety of length scales. For our model, we take $I$ to be the region within $d=1nm$ of the theoretical sharp interface. Note that if we take $I$ to be a straight line parallel to the contacts of the device and assume homogeneity in the parallel direction, we exactly recover the model in Ref. \refcite{BRG}.

We supplement the system \eqref{eq:basicelectrons} to \eqref{eq:basicpoisson} with the following boundary conditions. Dividing the boundary $\Gamma$ into  disjoint Dirichlet and Neumann parts $\Gamma=\Gamma_D\cup\Gamma_N$, we shall consider 
\begin{equation} \label{eq:boundaries}
\left\{ 
\begin{array}{llll}
V = V_D, & n = n_D, & p = p_D   & \quad\text{on}\quad \Gamma_D, \\
\nabla_\nu V = 0, & \nabla_\nu n = 0, & \nabla_\nu p = 0  & \quad\text{on}\quad \Gamma_N,
\end{array}
\right.
  \quad\text{and}\quad \nabla_\nu X = 0 \quad\text{on}\quad \Gamma \text{.}
\end{equation}
The boundary conditions for $\Gamma_D$ in Eq. \eqref{eq:boundaries} correspond to the conducting ends of the device (the left and right of Fig. \ref{fig:schematic}). The difference in the values for the potential $V$ at the two Dirichlet contacts correspond to the potential difference in the device (in Volts) whose offset does not affect the dynamics of the device because only the electric field $E=-\nabla V$ enters the other equations. The usual notation for a photovoltaic device follows the convention scheme of a diode, but in the working case current flows in the direction opposite to that of a usual diode. For this reason, a consistent applied potential is given at the left boundary of the device. We thus take the potential at the right side of the device to be zero, and take the potential at the left side to be the potential difference. Note that this notation has the convenient consequence that the average field in the device has the same sign as the potential difference.

For the boundary conditions of electrons and holes one has to take into account that the HOMO/LUMO levels of the device determine different energies for the electrons and holes implying that the concentrations of $n$ and $p$ at the metallic contacts can be very different. Thus, for majority carriers in their energetically favorable material ($p$ to the left of the interface and $n$ to the right as indicated in Fig. \ref{fig:schematic}), we take the boundary conditions given by Scott and Malliaras in Ref. \refcite{SM}. For the minority carriers ($n$ to the left and $p$ to the right), homogeneous Dirichlet boundary conditions constitute a very good approximation as the Scott-Malliaras boundary conditions give values approximately six orders of magnitude smaller than typical concentrations (see Table 1). % Explicit reference to \ref{table:parameters}).

On the insulated boundary $\Gamma_N$ we take homogeneous Neumann conditions for all of the variables because no particles should be able to pass into the insulator and it should be electrically neutral. For the excitons, we take homogeneous Neumann conditions on the whole boundary $\Gamma$ because the excitons are not be able to transmit from the organic material into the inorganic contacts.

The results of this paper is as follows: In Sec. \ref{sec:scaling} we will scale the physical parameters of the equations and detail the specific expressions for the mobilities, dissociation rates and recombination rates in an organic semiconductor material.  Next, in Sec. \ref{sec:asymptotics} we perform an asymptotic analysis of the 1-D equations for some typical parameters. The asymptotic analysis quantifies approximately how the exciton production rate yields an electric current and gives an expression for the shunt resistance of the device. Moreover, we present the numerical system used for the 2-D simulations in Sec. \ref{sec:numerics}. In Sec. \ref{sec:results} we plot the calculated concentrations in the device at three important points in the device operation - short circuit, optimal power point, and open circuit. We interpret these results in terms of their relation to the asymptotic results and furthermore show the current-voltage characteristic curves generated by both the numerical and asymptotic methods.

\section{Scaling and Models for Mobilities, Recombination, Dissociation} \label{sec:scaling}
We will scale \eqref{eq:basicelectrons} - \eqref{eq:basicpoisson}. The basic scaling introduces the following dimensionless quantities:
\begin{equation} \label{eq:basicscaling}
  V = \frac{k_b \mathrm{T}}{q} \tilde{V} = U_\mathrm{T} \tilde{V}, \qquad x = L \tilde{x}, \qquad (n,p,X) = N_r (\tilde{n},\tilde{p}, \tilde{X}), \qquad \mu = \mu_0 \tilde{\mu},
\end{equation}
where $U_\mathrm{T}$ is the thermal voltage, $L$ is a characteristic length (usually on the order of the device length), $N_r$ is a reference concentration, and $\mu_0$ is a reference mobility. We assume the Einstein's relation $D = U_\mathrm{T} \mu$. Although this might not be justified for some organic materials, it greatly simplifies the equations and as of yet a generally accepted alternative model has not been found. We further introduce the dimensionless quantities:
\begin{equation}
 \lambda^2 = \frac{\epsilon_0 U_\mathrm{T}}{q L^2 N_r}, \qquad T = \frac{L^2}{\mu_0 U_\mathrm{T}},
  \label{eq:lambdaTK} 
\end{equation}
where $T$ denotes a reference time scale. 
Note that the Debye length $\lambda$ is typically not a small parameter - see Table 1. %Explicit Reference because \ref{table:parameters} gives wrong number. 

For an OPV device under light we follow Ref. \refcite {BRG} and choose the Langevin recombination rate
$R_{np} = \frac{q(\mu_n + \mu_p) n p} {\epsilon_0 \epsilon_r} $, which rescales with the reference time $T$ as follows:
\begin{align*}
\frac{L^2}{N_r \mu_0 U_\mathrm{T}} R_{np} &= \frac{1}{\lambda^2}\frac{(\tilde{\mu}_n + \tilde{\mu}_p) \tilde{n} \tilde{p}}{\epsilon_r} = c_r \tilde{n} \tilde{p} \text{.}
\end{align*}
Setting $c_r' = c c_r$
this gives the following system:
\begin{align}
\lambda^2 \nabla \cdot (\epsilon_r \nabla \tilde{V}) &=  \tilde{n}- \tilde{p}-  \frac{h \chi_I}{L} \nabla_\nu \tilde{X} \label{eq:poissont} \\ 
\frac{\partial \tilde{n}}{\partial \tilde{t}} & = \nabla \cdot ( \mu_n \nabla \tilde{n} - \mu_n \tilde{n} \nabla (U+\tilde{V})) - c_r \tilde{n} \tilde{p} +  k_d T \tilde{X} \label{eq:electronst}\\
\frac{\partial \tilde{p}}{\partial \tilde{t}} &= \nabla \cdot ( \mu_p \nabla \tilde{p} + \mu_p \tilde{p} \nabla (U+\tilde{V})) - c_r \tilde{n} \tilde{p} +   k_d T \tilde{X} \label{eq:holest}\\
\frac{\partial \tilde{X}}{\partial \tilde{t}} &= \nabla \cdot (\mu_X \nabla \tilde{X}) + c_r' \tilde{n} \tilde{p} + GT - k_d T\tilde{X} - k_r T \tilde{X} \label{eq:excitonst} \text{.}
\end{align}

For the charge carrier mobilities $\mu_n$ and $\mu_p$ we shall apply the Poole-Frenkel model\cite{PF} postulating an exponential dependence on the square-root of the electric field $\tilde{E}=-\nabla \tilde{V}$:
\begin{equation}
\mu_n = \mu_n(0) e^{\gamma_n \sqrt{|\tilde{E}}|},\qquad
\mu_p = \mu_p(0) e^{\gamma_p \sqrt{|\tilde{E}}|}.
\end{equation}

Moreover, the exciton dissociation rate $k_d$ shall be given by Ref. \refcite{BRG} as a function of $M$, which is a scaled square-root of electric field $\tilde{E}$:
\begin{align*}
\left\{
\begin{array}{l} k_d^-(M) = \frac{2 k_d(0)}{M}\left( e^{M} \left(1-\frac{1}{M}\right) + \frac{1}{M}\right), \\[2mm]
k_d^+(M) = \frac{4 k_d(0)}{M^2} \left( 1 - e^{-M^2/4} \right),
\end{array}
\right.\qquad
M =\frac{1}{\lambda} \sqrt{\frac{|\tilde{E}|}{L^3 N_r \pi \epsilon_r}}, 
\end{align*}
where the $+$,$-$ superscripts represent positive and negative fields (with respect to the interface normal).

\subsection*{Device Parameter values}
The physical scaling parameters we take are:
\begin{align*}
 V &= U_\mathrm{T} \tilde{V} = .0258 \text{V}\, \tilde{V} \\
 E &= U_\mathrm{T} \tilde{E}/L = 2.58\times 10^5 \text{V/m}\, \tilde{E}  \\
 x &= L \tilde{x} = 100 \text{nm}\,\tilde{x}  \\
 n &= N_r \tilde{n} =  10^{20} \text{m$^{-3}$}\, \tilde{n}.  
\end{align*}

The values for the various constants are taken from Ref. \refcite{BRG} (except $G$ which is converted to a spatial density instead of an area density). Table 1 collects all the used dimensionless parameters.
%Manual insertion of '1' because \ref{table:parameters} gives '2' regardless of the placement of the \label{table:parameters} within the following table structure.

\begin{table}[htbp]
\tbl{Values for physical parameters.}{
\begin{tabular}{c c c c c} \toprule
$(h/L) = .01$  & & $(d/L) = .01$ & & $\epsilon_r = 4$\\
$\lambda^2 \approx 1.43$  & & $T = .00386$ && $E_{SC} = 13$  \\ 
$G T \approx 16990$ && $c_r=c_r' \approx .6987$ & & $U(x_r) - U(x_l) \approx 12$ \\ \colrule
$k_{d,out}T \approx 1$ & & $k_{r,in}T \approx 3.86$  & & $k_{r,out}T \approx  3864$ \\
$k_{d,in}(0)T \approx 386$ & & $k_{d,in}(+E_{SC})T \approx 178$ & & $k_{d,in}(-E_{SC})T \approx 2763$ \\ \colrule
$\mu_0 = 10^{-10}$ & & $(\mu_1 / \mu_0) \approx  .01$ \\
$\mu_n(0) = 3$ & & $\gamma_n = .788$ & & $\mu_n(E_{SC}) \approx 53.3$\\
$\mu_p(0) = 1$ & & $\gamma_p =.153$ & & $\mu_p(E_{SC}) \approx 1.75$\\ \colrule
$n(x_L),p(x_0) \approx .04$ & & $n(x_0),p(x_L) \approx 4 \times 10^{-7}$ \\ \botrule
\end{tabular}}
\begin{tabnote}
$E_{SC} = 3.5 \times 10^6$, corresponding to the short circuit potential difference of $-.5V$. The expressions $k_{out}$ and $k_{in}$ refer to the rate-values in- and outside of the interface region.
\end{tabnote}
\label{table:parameters}
\end{table}

The values of $\mu_X$ and $k_{d,out}$ are not easy to calculate physically and no consensus seems to exist in the literature for their values, but we take the given values based on reasonable estimates given in the sources listed above.

\section{Steady-State Equations}
For usual device operation, a solar cell will be producing current steadily for time-scales on the order of hours. Any transient behavior occurs over such short time-scales (i.e. microseconds\cite{HMG}) that we neglect them for modeling the long-term behavior and efficiency of the device. 

Therefore, for the remainder of the paper, we will consider only the steady-state equations where we drop all the tildes and absorb the time-scaling $T$ into the rates:
\begin{align} \label{eq:SteadyState}
\begin{cases}
 \lambda^2 \nabla \cdot (\epsilon_r \nabla V) &=  n- p-  \frac{h \chi_I}{L} \nabla_\nu X \\ 
-\nabla \cdot ( \mu_n \nabla n - \mu_n n \nabla (U+V)) &= - c_r n p +  k_d X \\
-\nabla \cdot ( \mu_p \nabla p + \mu_p p \nabla (U+V)) &= - c_r n p +   k_d X \\
-\nabla \cdot (\mu_X \nabla X)  &= c_r' n p + G - k_d X - k_r X \text{.} 
\end{cases}
\end{align}

We expect that this system of equations will have a unique steady-state solution. If we insist on the physically reasonable requirement that $k_d$ and $\mu$ are smooth and bounded from above (corresponding to physical device saturation), then the first three equations fit very nearly into the standard semiconductor framework (see, for example, Ref. \refcite{M}). The only remaining difficulty is the exciton equation, but given $n,p \in H^1 \cap L^\infty$ as in the usual case, we see that for smooth $\mu_X$ we have $X \in W^{1,\infty}$ and the exciton terms in the first three equations should not pose any difficulty. Proofs of the existence and uniqueness and the corresponding results for the time-dependent are under current investigation and subject to an upcoming paper.

Concerning the steady-state solutions, we emphasize that the light generation term $G > 0$ implies a non-zero right-hand-side in the n and p equations of \eqref{eq:SteadyState} and thus non-constant drift-diffusion fluxes of $n$ and $p$ on the left-hand-side. Thus, we expect a steady flux of electrons and holes away from the interface with nearly negligible bimolecular recombination (due to the work-function considerations discussed in the introduction).
Even for $G = 0$, we will only recover constant stationary drift-diffusion fluxes in such particular cases such as $\mu_X = 0$ and $c_r = \frac{c'_r k_d}{k_d+k_r}$, which is not realistic given the considered parameters. As a consequence, the \emph{model system \eqref{eq:SteadyState}  is not designed
to accurately predict the behavior of the device in the dark} $G = 0$. 

In the following section, we shall investigate further the steady state solutions
of \eqref{eq:SteadyState}  for $G > 0$ in a 1-D setting.

\section{1-D Stationary State Asymptotics} \label{sec:asymptotics}

\subsection{Large Applied Field Approximation} \label{sec:asymptoticsinV}
We begin by considering the steady state equations  in 1-D:
\begin{equation} \label{eq:1DSteadyState}
\begin{cases}
 \lambda^2 \epsilon_r V_{xx} &=  n- p-  \frac{h \chi_I}{L} X_x \\ 
- ( \mu_n n_x - \mu_n n (U+V))_x &= - c_r n p +  k_d X \\
- ( \mu_p p_x + \mu_p p (U+V))_x &= - c_r n p +   k_d X \\
-\mu_X X_{xx}  &= c_r' n p + G - k_d X - k_r X  \text{,}
\end{cases}
\end{equation}
where we have assumed that the $\epsilon_r$ and $\mu_X$ are spatially homogeneous.

In reverse bias, which constitutes the operating state of an OPV bilayer device, a negative 
voltage $V_{\mathrm{diff}}<0$ is given at the left boundary of the device, i.e. $V(x_0)=V_{\mathrm{diff}}$, while zero voltage is given at the right boundary of the device, $V(x_L)=0$. Note that with this notation, $V_{\mathrm{diff}}$ is the potential difference in the device. For working photovoltaic cells, this potential difference is the sum of the potential applied to the device and a built-in potential from the metallic contacts.

With a constant relative permittivity $\epsilon_r$ we introduce the Debye length  $\lambda_D^2:=\lambda^2 \eps_r\approx 5$, 
and rescale the potential according to this potential difference, i.e. $V\to |V_{\mathrm{diff}}| V$ which rewrites the Poisson equation from \eqref{eq:1DSteadyState} as
\begin{equation}
 V_{xx} =  \eps\left(n- p-  \frac{h \chi_I}{L} X_x\right), \qquad \eps:=\frac{1}{\lambda_D^2 |V_{\mathrm{diff}}|},
 \qquad  V(x_0)=-1, \ \  V(x_L)=0.
 \label{eq:poissonrescaled} 
\end{equation}
Since in a working device a typical short circuit voltage is about $-0.5$ Volts,\cite{BRG} which is 
$|V_{\mathrm{diff}}|\approx 19.3$ in our units, we find that $\eps\approx 0.01$. 

Hence (assuming for the moment that $\eps\frac{h}{L} X_x \ll 1 $,  where $\frac{h}{L}$ is the fraction of 
interface width to device width) we expect that the electric potential $V$ 
is in zeroth $\eps$-order dominated by the potential difference and, 
thus that the electric field is constant in the zeroth order of $\eps$, i.e. 
\begin{equation}
V^0(x) = V_{\mathrm{diff}} - E^0 (x-x_0), \qquad   E^0:= \frac{V_{\mathrm{diff}}}{x_L-x_0}.  \label{Vzero}
\end{equation}

We shall see in the following that the approximation \eqref{Vzero} remains consistent with
the assumption $\eps\frac{h}{L} X_x \ll 1 $ after being carried over to the equations for the charge carriers and the excitons .

However, first we quote from Ref. \refcite{SM} the typical order of magnitudes of the boundary values of electron- and hole- densities in the working state of the device:
\begin{equation*}
n(x_0) \ll \eps,\ \ n(x_L)=O(\eps), \qquad p(x_0) =O(\eps), \ \ p(x_L)\ll \eps \text{.}
\end{equation*}
Next, we remark that with the zeroth order approximation $E^0$ also the mobilities $\mu_n$ and $\mu_p$ are constant in zeroth. For the short circuit value $|E^0|=13$, we calculate that 
$\mu_n \approx 50$ and $\mu_p\approx 2$ and thus 
$1 / (\mu_n \mu_p) \approx 0.01 = \eps$. 
Thus, by rescaling 
\begin{equation}
\mu_n n \to n \quad\text{and}\quad \mu_p p\to p,
\label{eq:rescalenp}
\end{equation}
we obtain with $c_r=O(1)$
the following rescaled charge carrier equations 
\begin{align}
 -(n_x - n (V+U)_x)_x &= -O(\eps)\, n p + k_d X,\quad
n_0:=n(x_0) \ll \eps,\ n_L:=n(x_L)=O(1), \label{eq:electronsrescaled} \\
 -(p_x + p (V + U)_x)_x &= -O(\eps)\, n  p + k_d X,\quad
p_0=:p(x_0) =O(1), \ p_L:=p(x_L)\ll \eps.  \label{eq:holesrescaled}
\end{align}

Moreover, the rescaled equation for the excitons is
\begin{align}
 -\mu_X X_{xx} &= O(\eps)\, n p + G - (k_d + k_r) X, \qquad
 X_x(0)=X_x(x_L)=0,   \label{eq:excitonsrescaled}
\end{align}
where $\mu_X\ll 1$ is small, typically of order $O(\eps)$. 

Neglecting the quadratic $O(\eps)np$ term, we can solve this equation explicitly in terms of hyperbolic sines and cosines. This allows us to check the consistency the approximation \eqref{Vzero} and the rescaling \eqref{eq:rescalenp} with the underlying assumption $\eps\frac{h}{L} X_x \ll 1$: With $k_{r,out} \approx 4\cdot10^3$ and $k_{d,out} \approx 1$ outside the interface and with $k_{r,in} \approx 1$ and $k_{d,in} \approx 3\cdot10^3$ (for typical $E^0<0$) inside the interface we have:
\begin{align}
 \|X^0_x\|_{L^\infty([x_0,x_L])} \approx \frac{G}{\sqrt{\mu_X}}\left(\frac{1}{\sqrt{k_{d,in}}} - \frac{1}{\sqrt{k_{r,out}}}\right)
\end{align}
and thus $\eps\frac{h}{L} \|X^0_x\|_{L^\infty([x_0,x_L])} \approx O(\eps^{3/4})$ if $\frac{h}{L}=O(\eps)$ and $\mu_X=O(\eps)$.
 
\subsubsection{Zeroth Order Charge Carriers}
Now, we proceed to study the zeroth order solutions $n^0$ and $p^0$ of \eqref{eq:electronsrescaled} and \eqref{eq:holesrescaled}. Introducing the zeroth order fluxes 
\begin{align*}
 J^0_n := n^0_x - n^0 \varphi_x, \qquad\
 J^0_p := -(p^0_x + p^0 \varphi_x),\qquad\text{with}\quad \varphi=V^0+U,
\end{align*}
and neglecting the $O(\eps)$ recombination term, the equations \eqref{eq:electronsrescaled} and \eqref{eq:holesrescaled} integrate immediately as 
\begin{align}
   J^0_n(x) = J^0_{n0} - &F(x), \qquad   J^0_p(x) = J^0_{p0} + F(x), \label{eq:ChangeInFlux} \\ 
   &F(x):=\int_{x_0}^x{k_d(y) X^0(y) dy} \nonumber \text{.}
\end{align}

Utilizing the Slotboom variables $J^0_n = e^{\varphi}(n^0 e^{-\varphi})_x$ and $J^0_p = -e^{-\varphi}(p^0 e^{\varphi})_x$, we can explicitly solve for  $n^0$ and $p^0$:
\begin{align}
n^0(x) &= n_0 e^{\varphi(x)-\varphi_0} + J^0_{n0}\int_{x_0}^x{ e^{\varphi(x)-\varphi(y)}dy} - \int_{x_0}^x{F(y)e^{\varphi(x)-\varphi(y)}dy} \label{eq:ngeneral}\\
&=n_0 e^{\varphi(x)-\varphi_0} + J^0_{n0}\Phi_n(x) -\mathcal{F}_n(x),\nonumber\\
p^0(x) &= p_0 e^{\varphi_0-\varphi(x)} - J^0_{p0}\int_{x_0}^x{e^{\varphi(y)-\varphi(x)}dy} - \int_{x_0}^x{F(y)e^{\varphi(y)-\varphi(x)}dy} \label{eq:pgeneral}\text{,}\\
&= p_0 e^{\varphi_0-\varphi(x)} - J^0_{p0} \Phi_p(x)-\mathcal{F}_p(x),\nonumber
\end{align}
where $\varphi_0:=\varphi(x_0)$ and we define 
\begin{align*}
 \Phi_n(x) &:= \int_{x_0}^x{e^{\varphi(x)-\varphi(y)}dy}, &\qquad \Phi_p(x) &:= \int_{x_0}^x{e^{\varphi(y)-\varphi(x)}dy}, \\
\mathcal{F}_n(x) &:= \int_{x_0}^{x}{F(y)e^{\varphi(x)-\varphi(y)}dy}, &\qquad \mathcal{F}_p(x) &:= \int_{x_0}^{x}{F(y)e^{\varphi(y)-\varphi(x)}dy}.
\end{align*}

Next, we can determine the parameters $J^0_{n0}, J^0_{p0}$ by evaluating the boundary conditions $n(x_L)=n_L$ and $p(x_L)=p_L$. Upon rearrangement, we have with $\varphi_L:=\varphi(x_L)$:
\begin{align}
J^0_{n0} &= \frac{n_L - n_0 e^{\varphi_L-\varphi_0}  + \mathcal{F}_n(x_L)}{\Phi_n(x_L)},\qquad
J^0_{p0} &= \frac{-p_L + p_0 e^{\varphi_0-\varphi_L}  - \mathcal{F}_p(x_L)}{\Phi_p(x_L)},\label{eq:JnpExplicit}
\end{align}
and obtain explicit formulas for $n^0$ and $p^0$:
\begin{align}
n^0(x) &= n_0 \left(e^{\varphi(x)-\varphi_0} -  e^{\varphi_L-\varphi_0}\frac{\Phi_n(x)}{\Phi_n(x_L)}\right) + n_L\frac{\Phi_n(x)}{\Phi_n(x_L)} + \mathcal{F}_n(x_L) \frac{\Phi_n(x)}{\Phi_n(x_L)} - \mathcal{F}_n(x) \label{eq:nexplicit}\\
p^0(x) &= p_0 \left(e^{\varphi_0-\varphi(x)}- e^{\varphi_0-\varphi_L}\frac{\Phi_p(x)}{\Phi_p(x_L)}\right) + p_L   \frac{\Phi_p(x)}{\Phi_p(x_L)} + \mathcal{F}_p(x_L)\frac{\Phi_p(x)}{\Phi_p(x_L)} - \mathcal{F}_p(x) \label{eq:pexplicit}
\end{align}
which satisfy both boundary conditions since $\Phi_{n/p}(x_0)= \mathcal{F}_{n/p}(x_0)=0$.

These equations, although explicit, do not yet intuitively present the leading order contributions. 
However, using the fact that $U$ is constant outside of the interface, one can explicitly calculate the integrals $\Phi_n$ and $\Phi_p$ (see \ref{ap:phi}).

In the following, we use these formulas to identify the leading contributions to $J^0_{n0}$ and $J^0_{p0}$  and thus $n^0$ and $p^0$. More precisely, we introduce the parameters 
\begin{equation}
\frac{1}{\delta} := e^{\varphi(x_l)-\varphi_0},\qquad  \frac{1}{\eta}:=e^{\varphi(x_r)-\varphi(x_l)}  \label{eq:deltaeta}
\end{equation}
and observe that $e^{\varphi_L-\varphi(x_r)}=\delta^{-2}$ for the considered device geometry since we have $2(x_l-x_0)= x_L-x_r$ (see Fig \ref{fig:schematic}). With $\varphi_x=-E^0+U_x$ and $E^0<0$ as in the working device, these parameters are small: $\delta \approx 10^{-3}$ and $\eta \approx 10^{-6}$ for $E^0=E_{SC}$ with the numerical values given in Table 1. %Explicit reference to \ref{table:parameters}. 

However, using these parameters directly in \eqref{eq:nexplicit} and \eqref{eq:pexplicit} does not directly yield insights into the behavior of $n^0$ and $p^0$. Because $\Phi_{n/p}(x)$ can change over many orders of magnitude, the behavior of $n^0$ and $p^0$ depends highly on competing exponential terms, none of which can be easily eliminated. One possible way to proceed splits the domain $[x_0,x_L]=[x_0,x_l]\cup(x_l,x_r)\cup[x_r,x_L]$ into three areas: left of the interface, the interface, and right of the interface. For each of these areas, we would obtain formulas for $n^0$ and $p^0$ of the form of \eqref{eq:ngeneral} and \eqref{eq:pgeneral} in terms of the values $n^0(x_l), n^0(x_r)$ and $p^0(x_l), p^0(x_r)$ and the boundary terms $n_0$, $n_L$ and $p_0$, $p_L$. Then, imposing continuity of $n^0$ and $p^0$ and continuity of the fluxes $J^0_{n0}$ and $J^0_{p0}$ at $x=x_l$ and $x=x_r$ would allow us to determine the values $n^0(x_l), n^0(x_r)$ and $p^0(x_l), p^0(x_r)$.

However, a simpler way to proceed considers the zeroth order currents, which are sufficient 
to understand the produced electric current. Using the explicit formulas for $\Phi_n$ and $\Phi_p$ 
in terms of $\delta$ and $\eta$ given in \ref{ap:phi} we have from \eqref{eq:JnpExplicit}:
\begin{align}
J^0_{n0} &= \frac{n_L - \frac{n_0}{\eta \delta^3}  + \mathcal{F}_n(x_L)}{\frac{1}{E^0} \frac{1}{\eta \delta^2}\left(1 - \frac{1}{\delta}\right) + \frac{1}{E^0-U_x(\theta)} \frac{1}{\delta^2} \left(1 - \frac{1}{\eta}\right) + \frac{1}{E^0} \left(1 - \frac{1}{\delta^2}\right)} \text{,} \label{eq:JnParameters} \\[1mm]
J^0_{p0} &= \frac{-p_L + p_0 \eta \delta^3  - \mathcal{F}_p(x_L)}{-\frac{1}{E^0} \eta \delta^2 \left(1 - \delta \right) - \frac{1}{E^0-U_x(\theta)} \delta^2\left(1 - \eta\right) - \frac{1}{E^0} \left(1 -\delta^2\right)} \text{,} \label{eq:JpParameters}
\end{align}
where $U_x(\theta)$ denotes a mean-value of $U_x$ within the interface.

Because $J_n + J_p = const$, the sum $J^0_{n0}+J^0_{p0}$ can be used to calculate the total current in the device (as predicted by the asymptotics). We shall plot and discuss the predicted relationship between the current and the applied field as asymptotic IV-curve
in Sec. \ref{sec:results} (see Fig. \ref{fig:AsymptoticIV}). Note that for the given plots we take $U$ to be piecewise linear so that $U_x(\theta)$ is explicitly defined.

\subsubsection{Short-Circuit Current} \label{sec:AsymptoticCurrent}
As mentioned earlier, the sign of $E^0$ determines if the parameters $\eta, \delta$ are large or small. For the short circuit current, we have $E^0 < 0$ and $\eta,\delta \ll 0$. In order to determine the lowest order terms of these equations, we must also calculate a more precise form of $\mathcal{F}_{n/p}(x_L)$ and thus its order. At first we remark that $F(x)=\int_{x_0}^x{k_d(y) X^0(y) dy}=O(10^2)$
 (since $k_{d,in} X^0 = O(G)=O(10^4)$ on the interface with width $10^{-2}$ and $k_{d,out} X^0 = O(10)$ outside the interface with device length $x_L-x_0 = O(1)$). Thus, one can verify with calculations similar to those below and in \ref{ap:phi} that the exponentials of $\varphi(x)$  will determine the leading contributions to the integrals $\mathcal{F}_{n/p}(x_L)$. More precisely, 
since $\varphi(x)$ is strictly increasing in the $E^0<0$ regime (recall that $\varphi_x=-E^0+U_x$) the leading order contributions derive from 
\begin{align*}
 \mathcal{F}_n(x_L) \approx \int_{x_0}^{x_l}{F(y)e^{\varphi_L-\varphi(y)}dy}, \qquad
 \mathcal{F}_p(x_L) \approx  \int_{x_r}^{x_L}{F(y)e^{\varphi(y)-\varphi_L}dy} \text{.}
\end{align*}
Applying the definitions of Eq. \eqref{eq:deltaeta}, integration by parts and the mean value theorem yields
\begin{align*}
\int_{x_0}^{x_l}{F(y)e^{\varphi_L-\varphi(y)}dy} 
&=\frac{1}{\eta \delta^2}  \int_{x_0}^{x_l}{F(y)e^{-E^0(x_l-y)}dy}  \\
 &= \frac{1}{\eta \delta^2}\left[ F(y) \frac{e^{-E^0(x_l-y)}}{E^0}\right]_{x_0}^{x_l} - \frac{1}{\eta \delta^2}\int_{x_0}^{x_l}{\frac{F'(y)}{E^0}e^{-E^0(x_l-y)}dy}  \\
 &= \frac{1}{\eta \delta^2} \frac{F(x_l)}{E^0}-\frac{1}{\eta \delta^2}\frac{F'(\theta_n)}{(E^0)^2}\int_{x_0}^{x_l}{e^{-E^0(x_l-y)}dy}, 
\end{align*}
for a mean value $\theta_n\in(0,x_l)$. Together with similar calculations for $\mathcal{F}_p(x_L)$ and $\theta_p\in(x_r,x_L )$, we obtain 
\begin{align*}
 \mathcal{F}_n(x_L)&\approx \frac{1}{\eta \delta^2} \frac{F(x_l)}{E^0} -\frac{k_{d,out}\, X^0(\theta_n)}{(E^0)^2}\frac{1}{\eta \delta^2}\left(1 - \frac{1}{\delta}\right)  + O(\eta^{-1} \delta^{-2})\\[1mm]
 \mathcal{F}_p(x_L) &\approx -\frac{k_{d,out}\, X^0(\theta_p)}{(E^0)^2}\left(1 - \delta^2\right) - \frac{F(x_L)}{E^0} + \frac{F(x_r)}{E^0} \delta^2 + O(\delta^2)
\end{align*}

 Taking only the dominating terms $O(\eta^{-1}\delta^{-3})$ in \eqref{eq:JnParameters} and ${O(1)}$ in \eqref{eq:JpParameters} yields:
\begin{align*}
 J^0_{n0} &\approx E^0 n_0 - \frac{k_{d,out} X^0(\theta_n)}{E^0} \qquad J^0_{p0} \approx E^0 p_L - \frac{k_{d,out} X^0(\theta_p)}{E^0} - F(x_L)\\[1mm]
 J^0_{tot} &\approx E^0 (n_0 + p_L) - \frac{k_{d,out}}{E^0}\big(X^0(\theta_n)+X^0(\theta_p)\big) - \int_{x_0}^{x_L}{k_d X^0(y) dy}\text{.}
\end{align*}

The term $E^0(n_0+p_L)$ gives a current which is proportional to the field, indicating that it represents a resistor. For a usual working device this term will be small ($n_0+p_L \approx 10^{-6}$ in Sec. \ref{sec:scaling}). 
Because the term does not depend on the length of the device, we interpret it as a shunt resistance (parallel to the device). We observe from Eqs. \eqref{eq:JnParameters} and \eqref{eq:JpParameters} that we have another current term with the form of a shunt resistance: $\eta \delta^3 E^0 (n_L + p_0)$. This factor is negligible in the short circuit case, but becomes large as we move to positive fields, where $\delta^3 > 1/\eta$. We investigate the total shunt resistance further in Sec. \ref{sec:results}. 

The second two terms represent the interaction of the excitons on the system. For usual device parameters, $k_{d,out} \ll k_{d,in}$ (for the parameters given in Sec. \ref{sec:scaling}, we have $k_{d,out} = 10^{-3} k_{d,in}$). Since $X^0=O(\frac{G}{k_{d,in}})$ on the interface and $X^0 = O(\frac{G}{k_{r,out}})$ outside the interface, we can neglect the contributions from outside the interface (recall that $\theta_n\in(0,x_l)$ and $\theta_p\in(x_r,x_L )$) compared to the contribution from within the interface.

Thus the most important contribution to the current for negative fields is:
\begin{equation} \label{eq:JApprox}
 J^0_{\mathrm{approx}} = - \int_{x_l}^{x_r}{k_{d,in} X^0(y) dy} \text{.}
\end{equation}
This approximation works very well in the short circuit case (as examined in Sec. \ref{sec:results}, in the discussion preceding Fig. \ref{fig:SC}), where it replicates $J^0_{tot}$ within 2\% and is actually closer to the numerically calculated $J$.

Note that $J^0_{\mathrm{approx}}$ comes entirely from the $J^0_{p0}$ term. Because we have evaluated $J^0_n$ and $J^0_p$ at the point $x_0$, we expect that the hole current will be dominant. However, using Eq. \eqref{eq:ChangeInFlux} we can calculate $J^0_{nL}$ and $J^0_{pL}$, and that the only change from the currents at $x_0$ is that the $-F(x_L)$ term appears in the electron current. Thus, as expected, in the region where the electrons are favored, the primary contribution to the current comes from the electrons and the total current is conserved.

\subsubsection{First Order Terms} \label{sec:firstorder}
Since we have explicit forms for $n^0$, $p^0$, and $X^0$, we can integrate twice to calculate $V^1$ (with an additional linear term to account for the two boundary values of $V$). Normally we don't plot $V$ since it is generally dominated by its boundary terms. However, $E^1 = - V^1_x$ can also be calculated explicitly (by a single integration) and we include the first-order field in the plots in Sec. \ref{sec:results}. The explicit form is not very illuminating and thus not written here. Note that we need to add a constant value to the field to insure that its integral gives the correct potential difference in the device (corresponding to the slope of the aforementioned linear term).

In theory we could use the second order form for the potential to calculate the next order $n^1$ and $p^1$ solutions. In fact, this is more or less the essence of the iteration scheme outlined in Sec. \ref{sec:numerics}. However, without a constant field, it is no longer possible to explicitly integrate the continuity equations (especially given the sub-exponential forms for the electron and hole mobilities). It would be possible to do these calculations numerically, but this would simply become a simplified version of our numerical iteration scheme, and thus we do not pursue this further in the asymptotic case. Instead, we will present another method for asymptotic calculations.

\subsection{Unipolar approximation} \label{sec:unipolar}

The discussion of the previous section, in particular Eq. \eqref{eq:JApprox}, can be summarized by the statement that the electron- and hole- fluxes $J^0_n$ and $J^0_p$ are approximately constant outside the interface but feature a strong variation over the interface with a magnitude of $F(x_r)-F(x_l)=\int_{x_l}^{x_r}{k_{d,in} X^0(y) dy}$. 

As a consequence, using the Eqs. \eqref{eq:nexplicit} and \eqref{eq:pexplicit} and similar calculations as in Sec. \ref{sec:AsymptoticCurrent} and \ref{ap:phi}, it follows that the zero-order electron- and hole-densities $n^0$ and $p^0$ vary also strongly over the interface. In fact, we obtain for the hole density $p^0$:  
\begin{align*}
 p^0(x_l) &= \frac{F(x_L) - F(x_l)}{-E^0} + O(\delta) + O(p_L)\\
 p^0(x_r) &= \frac{F(x_L) - F(x_r)}{-E^0+U_x(\theta_{in})} + \frac{k_{d,in} X^0(\theta_{in})}{(E^0-U_x(\theta_{in}))^2} + \frac{k_{d,out} X^0(\theta_{out})}{E^0(-E^0+U_x(\theta_{in}))} \\
 & \quad + O(\eta,\delta^2)+ O(p_L),
\end{align*}
where $\theta_{in} \in(x_l,x_r)$ is an mean value within the interface, $\theta_{out} \in (x_r,x_L)$ is a mean value outside of the interface and $\eta$, $\delta$ and $p_L$ are (as mentioned above) small in the working device. Then, since $U_x(\theta) \gg 1$ for $\theta\in(x_l,x_r)$ and because the dominant part of the integral $F(x)$ comes from the interface region, we see that 
$$
\frac{p^0(x_l)}{p^0(x_r)}= O(10^3), \qquad \text{in the working case }
E^0<0, 
$$
as a consequence of the work-function gap over the interface.
A similar result can be obtained for $n^0$ giving $n^0(x_r) \gg n^0(x_l)$. 

Note that in a physical device one might expect this ratio to be even larger. If we narrow the interface while keeping the work-function difference $(U(x_r)-U(x_l)$ and the generation term given in \eqref{eq:JApprox} constant, we tend to increase the value of $U_x(\theta_{in})$. Thus in some sense the smallness of $p^0(x_r)$ increases naturally from reducing the size of the interface.
 
The previous estimates quantify the observation that a bilayer device is operated such that in one material, $n \ll 1$ (here to the right of the interface) and in the other material $p \ll 1$ (here to the right of the interface) because of the nature of the HOMO/LUMO levels of the different polymers.  We can thus ask for an asymptotic simplification of system \eqref{eq:1DSteadyState} in the bulk-material away from the interface. 

In order to proceed with such an approximation, we observe furthermore that $\frac{p^0(x_l)}{k_{d,out} X^0}\approx 10 - 20$, i.e. that the boundary value of the majority charger carrier $p^0(x_l)$ given at the left of the interface dominates over the generation of holes from excitons in the bulk material left of the interface. In fact the above factor would even be larger in more realistic exciton models compared to \eqref{eq:1DSteadyState}, which misses a term modeling a electrochemical trapping of excitons (due to stability of polaron pairs on the interface\cite{AHB}) yielding further decreased the exciton concentrations. 
We shall thus entirely neglect here the photogeneration term of charge carriers outside the interface ($k_{d,out} X^0$).

Altogether, neglecting the minority charge carriers and excitons outside the interface, we consider the following 
simplified models for the majority charge carriers:
\begin{align}
 -\lambda_D^2 V_{xx} = \lambda_D^2  E_x = p, \qquad\qquad
 p_x - p E = \frac{-J_{p0}}{\mu_p}, 
\label{eq:HolesAiry}
\end{align}
in the $p$-material at the left of the interface (with the Debye length $\lambda_D^2=\lambda^2 \epsilon_r$) and
\begin{align*}
-\lambda_D^2  V_{xx}= \lambda_D^2  E_x = -n, \qquad\qquad
 n_x + n E = \frac{J_{n0}}{\mu_n},
\end{align*}
for the $n$-material at the right of the interface. Here $J_{p0}$ and $J_{n0}$ are the constants arising from integration. 

For the rest of this section we shall focus only on the hole density $p$. The analysis of the electron equation is analogous. 
We remark that the Eqs. \eqref{eq:HolesAiry} constitute a generalization of the system considered in Appendix A of Ref. \refcite{Cetal}, where the zero-current case $J_{p0}=0$ was analyzed. For the present case with non-zero current we assume first that $\mu_p$ is constant and exploit then a first-integral 
of the equation to obtain
\begin{equation*}
E_x - \frac{1}{2} E^2 = \frac{-J_{p0}}{\lambda_D^2  \mu_p} x + C_p
\end{equation*}
where $C_p$ is another integration constant. Using the  substitution $x = -2u$ and $E= y_u / y$ (see e.g. Ref. \refcite{Cetal}), we obtain
\begin{equation*}
 y_{uu} = (\kappa u - 2 C_p) y
\end{equation*}
where $\kappa = \frac{-4 J_{p0}}{\lambda_D^2  \mu_p}$. With the further substitution $z = \sqrt[3]{\kappa}\left( u - \frac{2 C_p}{\kappa}\right)$, this becomes the familiar Airy Equation, $y_{zz} = yz$ and we can write down the formula for $y$ explicitly in terms of Airy Functions. Since $E = y_u/y$, it is more convenient to write the form of $V = -\log{(y)}$ (derived independently in Ref. \refcite{RPK}):
\begin{equation}
V = -\log{\left[ \alpha Ai\left[\sqrt[3]{\kappa}\left( -\frac{x}{2} - \frac{2 C_p}{\kappa}\right)\right] + \beta Bi\left[\sqrt[3]{\kappa}\left( -\frac{x}{2} - \frac{2 C_p}{\kappa}\right)\right]\right]}
\end{equation}
where we can express $\alpha$ and $\beta$ in terms of the integration constants and the boundary values. Note that for the $n$ material we take the substitution $x=2u$ and $J_n$ positive and everything else proceeds identically.

This asymptotic result gives an explicit form for the charge carriers in the bulk in 1-D. However, for bilayer systems away from the interface, the system generally approximates a 1-D system since the bulk material acts primarily as a charge-transport layer. Since we explicitly account for the self-consistent potential, this formulation works much better than our previous asymptotics when we have large carrier concentrations (especially for less negative potential differences). In theory we could combine these equations in the bulk with a numerical calculation of the interface, simplifying the computations. However, in practice, we can change the mesh (or grid) to account for the reduced complexity of the problem in the bulk (see Fig. \ref{fig:Mesh}) and the gain would not necessarily be significant.

The main disadvantage of this formulation is that it relies upon the boundary data of the system. Since it doesn't model the interface, it is thus incapable of predicting the current in the device based on the device parameters, and requires $J$ as a parameter. In fact, we need all of the boundary values $p_0, E_0, J_{p0}$ in order to calculate the field and the hole concentration in the $p$-region of the device, but generally only $p_0$ is given. The other values, $E_0$ and $J_{p0}$, both require calculating the solution of the whole problem or at least calculating a model of the interface and using asymptotic approximations of the boundary values $E_0$ and $J_{p0}$. Alternatively, we can take the values $p(x_l), E(x_l), J_{pl}$ and again calculate the values in the bulk. Either method gives extremely accurate results for the cases in which $J \ll 0$, which we show in Fig. \ref{fig:SC} - \ref{fig:OPP}.

\section{Numerical Scheme} \label{sec:numerics}
\subsubsection*{Time-discretization}
We present a numerical scheme based on a Gummel-type iteration,\cite{gummel} which is well-suited to deal with the  nonlinear dependence of the equations for $n$, $p$, and $X$ on the electric field ($E = -\nabla V$). At each Gummel-iteration step we calculate first the electric potential $V$ from the current state of $n,p,X$ and then use this new electric potential to update $n,p,X$. We can write this as follows after the $k^{th}$ step:
\begin{equation} \label{eq:numericalsystem}
\begin{split}
&-\lambda^2 \nabla \cdot (\epsilon_r \nabla V_{k+1}) =  p_{k}- n_{k} +  \frac{h}{L} \nabla X_{k} \\
&-\nabla \cdot ( \mu_n(E_{k+1}) \nabla n_{k+1}  - \mu_n(E_{k+1}) n_{k+1} \nabla (V_{k+1}+U)) + c_r n_{k+1}\, p_k  =  k_d(E_{k+1}) X_k \\
&-\nabla \cdot ( \mu_p(E_{k+1}) \nabla p_{k+1}  + \mu_p(E_{k+1}) p_{k+1} \nabla (V_{k+1}+U)) + c_r n_k\, p_{k+1}  =  k_d(E_{k+1}) X_k \\
&-\frac{\mu_1}{\mu_0} \nabla \cdot (\mu_X \nabla X_{k+1}) + (k_d(E_{k+1}) + k_r) X_{k+1} = c_r' n_k p_k + G \text{.}
\end{split}
\end{equation}

Note that we use a semi-implicit discretization in time. Specifically, we implicitly discretize the flux and mass terms in \eqref{eq:numericalsystem}, but use semi-implicit discretization for the recombination term $c_r n p$ and an explicit form for the dissociation term $k_d X$. Thus $n$, $p$, and $X$ can be effectively updated in parallel since they only require the previously updated values for $V$. Furthermore, each variable $n$, $p$, $X$ is treated completely implicitly within its equation, which greatly reduces numerical instabilities arising from the errors in the previous step.

The iteration continues until a preset $L^2$ error between the two latest steps is achieved. 
In practical numerical simulation, we have found this Gummel iteration to converge well except for high positive values for the applied potential difference (which constitutes the non-working case). This case, which is characterized by opposite convection terms coming from the work-function and the electric potential, leads to large values for $n$ and $p$ near the interface which prevent the Gummel iteration from converging. To bypass this difficulty, we add a damping parameter $\alpha$ which interpolates between old and new solutions. More precisely, denoting with $V'_k,n'_k,p'_k,X'_k$ the solutions of the above Eqs. \eqref{eq:numericalsystem}, the damped scheme is:
\begin{equation} \label{eq:numericaldamping}
\begin{split}
 V_{k+1} &= \alpha V'_{k+1} + (1 - \alpha) V_k, \qquad\qquad n_{k+1} = \alpha n'_{k+1} + (1 - \alpha) n_k, \\
 p_{k+1} &= \alpha p'_{k+1} + (1 - \alpha) p_k, \qquad\qquad X_{k+1} = \alpha X'_{k+1} + (1 - \alpha) X_k \text{.}
\end{split}
\end{equation}
In particular, note that since $V'_{k+1}$ and $V_k$ both satisfy linear boundary conditions, so does any convex combination of $V'_{k+1}$ and $V_k$. This is equally true for the boundary conditions for $n$,$p$, and $X$, which are fixed for a given potential $V_{k+1}$. 

Our numerical experiments in the case of positive applied potentials show that it is possible to optimize the convergence by tuning the damping parameter. In fact it turns out that using a small damping $\alpha \approx .01$ for the first few (e.g. three) iteration steps serves to correct the initial guess and prevents the concentrations from going negative due to badly chosen initial data. Further optimization on the damping, such as a potentially dynamic choice of $\alpha$, is a topic for further investigation.

\subsubsection*{Space-discretization}
Our choice of a space-discretization has to take into account the two major difficulties inherent in our system. First of all, 
the material interfaces will cause (as discussed in the previous section) very strong changes in the densities $n$, $p$, and $X$ over and on the interface. Secondly, depending on the electric field, the equations can be either convection or diffusion dominated. Furthermore, we want our scheme to easily generalize to complicated 2-D interfaces.

The Hybrid Discontinuous Galerkin (HDG) finite element method is well-suited to handle these challenges. The novel idea for the HDG method is to have separate degrees of freedom on the elements $u$ (an element in the usual sense) and on the facets (boundaries) of the elements $u_F$. We assume that our domain is decomposed into a mesh $\mathcal{T}$ consisting of triangles with a set of facets $\mathcal{F}$ consisting of the facets (edges) of these triangles. The space we use for each of the equations is $V' = \{(u,u_F): u \in H^1(\Omega)\cap H^2(\mathcal{T}), u_F \in L^2(\mathcal{F})\}$. Our solution space is  $(V,V_f,n,n_f,p,p_f,X,X_f) \in V = V' \times V' \times V' \times V'$ and we denote the approximate solutions to the four equations as $u_V, u_n, u_p, u_X$ with corresponding facet functions $u_{Vf}, u_{nf}, u_{pf}, u_{Xf}$.  The facet elements allow for upwind-type calculations for convection-dominated behavior, as well as allowing for large changes of behavior at material boundaries. The functions within the elements are approximated by polynomials of a given order. Both increasing the polynomial order and refining the mesh improve the properties of a solution after a given number of iterations at the price of larger computational costs. 

A solid introduction to HDG methods is Ref. \refcite{HDG}, specifically the first section on the general convection-diffusion case. The HDG method can more accurately be called a Hybridized symmetric interior penalty Discontinuous Galerkin method. The relationship to standard Discontinuous Galerkin methods runs quite deeply. The stabilization and symmetrization terms have close analogs in traditional DG methods, see for instance Refs. \refcite{CGL}, \refcite{CKS}, \refcite{ABCM}.

The specifics of the weak form of the equations are also given in Ref. \refcite{HDG}. Here we shall recall only a few specific features. It is well-known that numerical methods for convection-dominated systems need to apply the correct upwinding. For our elements, this consists of choosing the facet elements for inflow boundaries, and the bulk elements for outflow boundaries. This selection allows for calculation along the flow of the carriers in the usual upwinding sense.

All the numerical examples consider the two-dimensional computational domain $\Omega= [0,1.5]\times[0,0.2]$ (in the scaled variables from Sec. \ref{sec:scaling}). In order to compare with the one-dimensional asymptotic analysis above, we take the domain interfaces (see Fig. \ref{fig:schematic}) to be straight vertical lines located at $x_0=0$, $x_l=0.49$, $x_m=0.5$, $x_r=0.51$, and $x_L=1.5$ respectively. 

The mesh generation was accomplished using the Netgen program created by Sch\"oberl.\cite{netgen} In particular, the mesh generation automatically preserves the internal interfaces allowing all coefficients to be defined piecewise depending on the properties of the local material. The numerics were done using NGSolve, a solver wrapper for the Netgen program. We also included the inverse solver PARDISO to assist with the non-symmetric matrix calculations.\cite{pardiso1,pardiso2} See Ref. \refcite{MW} for more details about using the NGSolve program on convection-diffusion equations.

\section{Discussion of Numerical Examples} \label{sec:results}

In the following we present selected numerical examples of the steady-state device behavior
of the organic photovoltaic bilayer device plotted in Fig. \ref{fig:schematic}. In particular, we shall 
investigate the steady-states calculated from different applied potentials.  

\begin{figure}[htbp] 
\centering
\scalebox{0.3}{\includegraphics{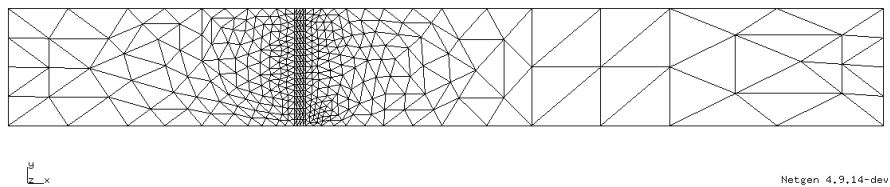}}
\caption{Mesh used for all of the plots in Sec. \ref{sec:results}} \label{fig:Mesh}
\end{figure}

Using the scheme outlined in Sec. \ref{sec:numerics} we use the damped-Gummel iteration 
given in \eqref{eq:numericalsystem} and \eqref{eq:numericaldamping} with $\alpha=.01$ for the first three steps (as a sort of preconditioning) and $\alpha = .6$ until convergence. All of the plots show a mesh consisting of 720 elements, greatly refined near the interface, and use polynomial interpolants of order 7 for the bulk elements - see Fig. \ref{fig:Mesh}. Refining the mesh and increasing the interpolation order both give better accuracy, in particular for the steep concentration changes over the interface (see, i.e. Figs. \ref{fig:SC}, \ref{fig:OPP} and \ref{fig:OC}). Note that 
the geometry of the test cases is taken to be homogeneous in the $y$-direction and we plot only a sample $x$-cross section which can be compared with the one-dimensional asymptotics. A fully 2-D numerical example is given at the end of the section.

\subsection{IV Curve} 
One of the most important characteristics of an OPV device is the effect of changing the applied voltage on the current in the device. We show this relationship by plotting the current as a function of potential difference by multiple simulations with changing boundary values. We pick values for $V_{\mathrm{diff}}$ and then use these values to calculate the boundary values for $n$ and $p$ according to Ref. \refcite{SM}. Due to the built-in potential of the metallic contacts at the OPV device, we have that zero  applied voltage leads to an approximate potential difference over the device of $V_{\mathrm{diff}} = -19.3$. 

\begin{figure}[htbp] 
\centering
\scalebox{.8}{\includegraphics{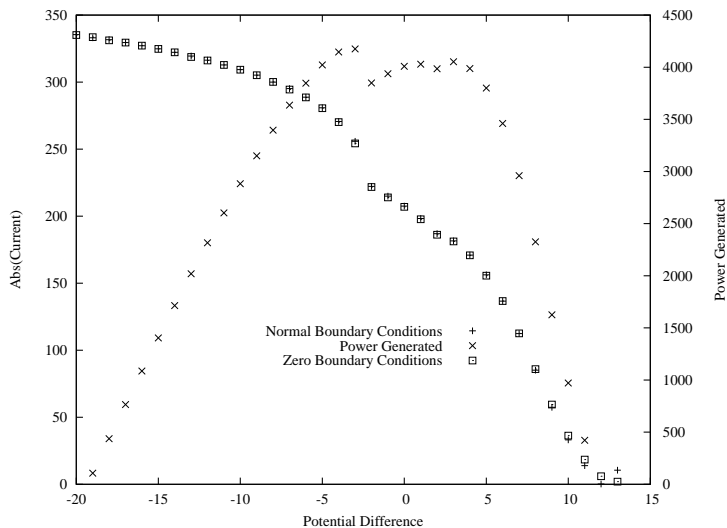}}
\caption{A comparison of the IV-curves obtained using the boundary conditions outlined above in Table 1 %explicit \ref{table:parameters} 
above and homogeneous Dirichlet conditions for the charge carriers and a plot of the generated power density for the usual boundary conditions. Note that the power corresponds to the areas of the rectangles with corners at the origin and the current at each applied potential.} \label{fig:IVBCs}
\end{figure}

Fig. \ref{fig:IVBCs} plots the IV-curve for our typical OPV device (marked with $+$) with respect to the potential difference. The zero applied voltage point is approximately at the left edge of the graph. This referential IV-curve compares very closely with a second IV-curve (marked with $\square$) computed after replacing the stated boundary conditions\cite{SM} by homogeneous Dirichlet boundary conditions for each carrier at both Dirichlet boundaries.

Both IV-curves show how the (absolute value of the) produced current decreases as the potential difference increases. In particular we observe a distinctive S-shape. This is an effect of the dependence of $k_{d,in}(E)$ on the self-consistent electric field $E$. The kink arises at the potential $V\approx-3$ for which the $E$-field touches zero from below at the interface. Note that the $E$-field touching zero happens for a potential $V$ smaller than zero because the holes on the left side of the device (i.e. the majority charges) increase self-consistently the negative applied $E$-field towards the interface. This argumentation can be verified by plotting the IV-curve obtained from a constant dissociation rate $k_{d,in} = k_{d,in}(E_{SC})$. A plot comparing this result to the asymptotic calculation is shown below in Fig. \ref{fig:IVkdconst}.

The power density generated by a device is given by $P = J (V_{\mathrm{int}}+V_{\mathrm{diff}})$. As shown in Fig. \ref{fig:IVBCs}, the power is zero for two specific characteristic points in the device: $V_{\mathrm{diff}} = -V_{\mathrm{int}}$ (short-circuit) and $J=0$ (open-circuit). A third characteristic point, the optimal power point, is defined by the voltage point with the maximal value of $P$. For potential differences outside of the regime shown in Fig. \ref{fig:IVBCs} the power generated is negative, indicating that the device requires power input to operate and is no longer in the working regime for a photovoltaic device. An upper bound for the maximum power is given by the product of the short circuit current and the open circuit voltage (these being the maximum current and maximum voltage for the device in the working regime). The \emph{fill factor} denotes the ratio of the maximum obtainable power and this upper bound.\cite{FF} Geometrically this is represented by the largest possible rectangle contained under the IV curve and the smallest rectangle containing the IV curve. For the plots given in Fig. \ref{fig:IVBCs} the maximal power point occurs at a potential difference of $V\approx-3$ with a fill factor of 39.5\%. This fill factor is in good agreement with other simplified organic devices.\cite{GEHW} In general, fill factors vary heavily depending on the device, and are notably improved for multijunction devices.

\subsection{Electron-/Hole- and Exciton- Densities at SC, OPP and OV}
In this section, we study the details of the charge density and exciton distributions at particular characteristic points of an IV curve: i) Short Circuit (SC), ii) Optimal Power Point (OPP) and iii) Open Circuit Voltage (OV). 
The short circuit current is the current at zero applied voltage, or a potential difference of $V_{\mathrm{diff}}=-19.3$. The optimal power point is the potential difference for which the power is maximized, which is $V_{\mathrm{diff}}=-3$ in Fig. \ref{fig:IVBCs}. The open circuit voltage is the point where $J=0$, which occurs for $V_{\mathrm{diff}}=12$.

For each of these points, we will plot the numerical solutions for $n$ (marked with $\square$) and $p$ (marked with $\times$) in the same chart. On the same figure, we also plot the zero-order asymptotics derived in Sec. \ref{sec:asymptoticsinV}. Note that we have rescaled $n$ and $p$ by $\mu_n$ and $\mu_p$ in the asymptotics section, and must appropriately scale them back for the plots. Moreover, we show how the unipolar asymptotics derived in Sec. \ref{sec:unipolar} reproduce the numerical solution in the bulk very well when complemented with the numerical data from the simulation of the interface (which is not part of the unipolar approximation). Specifically, we use the boundary data values $p(x_0)$, $E(x_0)$, and $J$ for the holes and $n(x_L)$, $E(x_L)$ and $J$ for the electrons. The unipolar approximation thus has the potential to reduce simulation costs for the dynamics of the bulk material. In particular we observe that in all three cases -- SC, OPP and OV -- the minority carriers ($n$ to the left of the interface and $p$ to the right of the interface) are strongly dominated by the majority carriers (as discussed above in Sec. \ref{sec:unipolar}).

Next, we plot the numerical solutions of the exciton concentration and the electric field.
The electric field is compared with the asymptotic approximation explained in Sec. \ref{sec:firstorder}. The unipolar asymptotics do not include exciton behavior to compare with the numerical results. Furthermore, we refrain from plotting the unipolar approximation to the electric field outside the interface (which is discontinuous over the interface since the field is calculated independently on the two sides of the device). However, good agreement similar to the charge carrier densities can be achieved. 

We also list the numerically calculated current given at each voltage. The current predicted by the zero-order asymptotics is the sum of $J_{nl}^0$ and $J_{pl}^0$ given in Eq. \eqref{eq:JnParameters} and \eqref{eq:JpParameters}. However, because this predicted current ignores recombination, it only gives a good approximation in the short circuit case.

\subsubsection*{Short circuit}
The case with the largest (negative) potential difference $V_{\mathrm{diff}}$ is the short circuit case (which has zero applied voltage). This is due to the built-in potential arising from the work functions of the anode and cathode metals. For the device parameters laid out in Sec. \ref{sec:scaling}, this occurs for $V_{\mathrm{diff}} \approx -19.3 $.

\begin{figure}[htbp] 
\centering
\scalebox{0.8}{\includegraphics{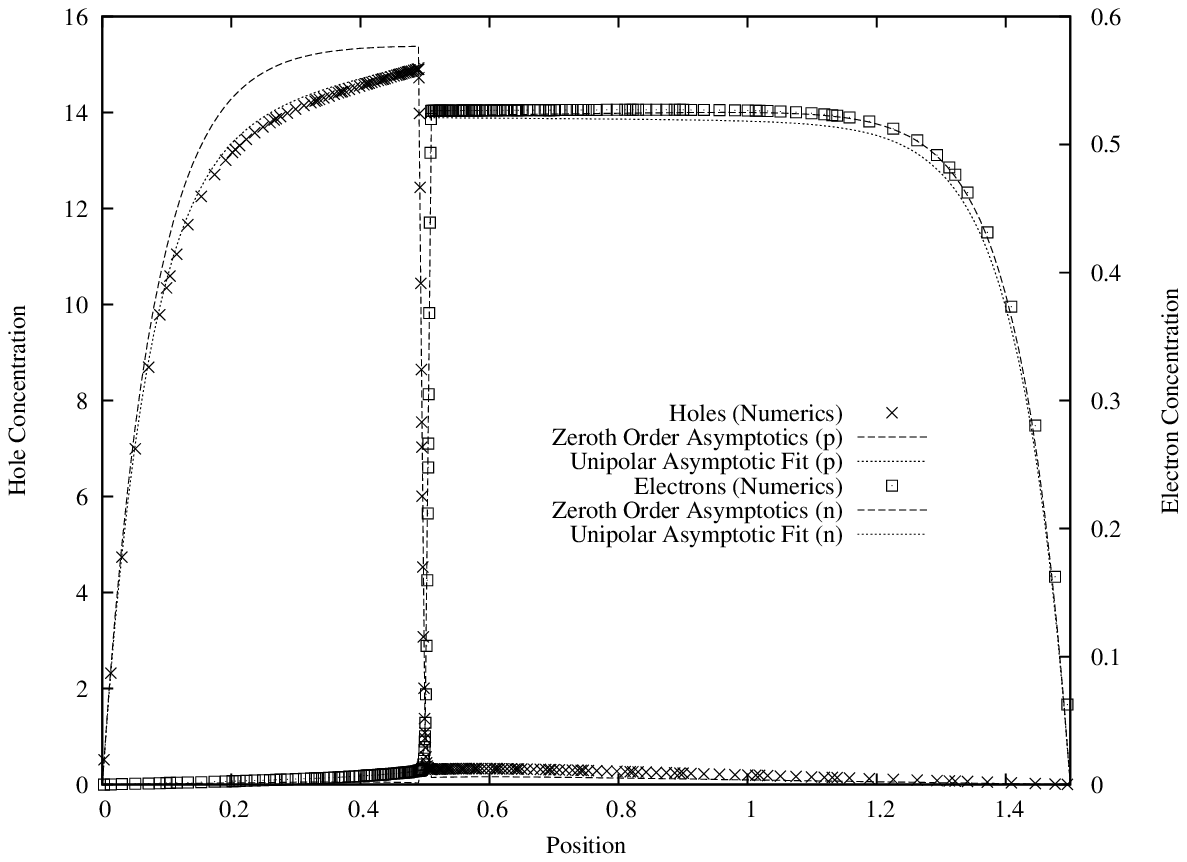}}
\scalebox{0.5}{ \includegraphics{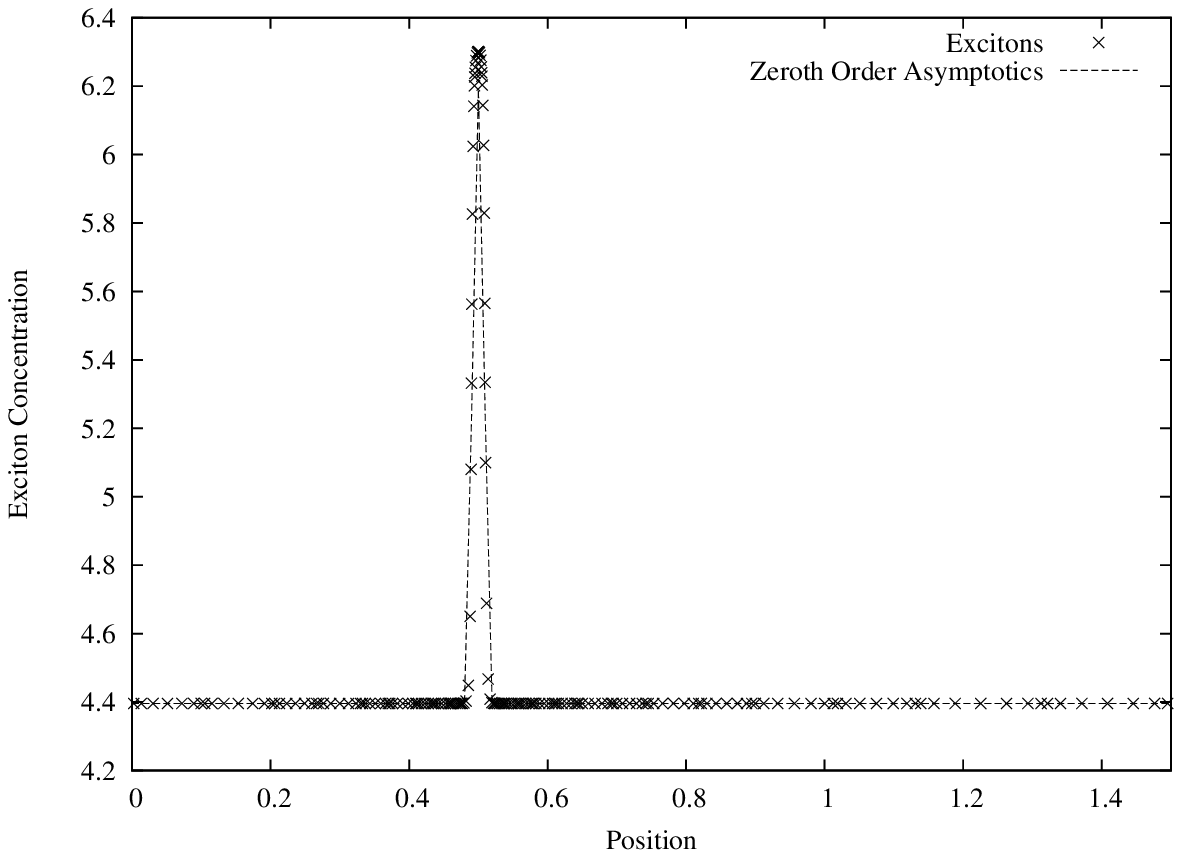} \includegraphics{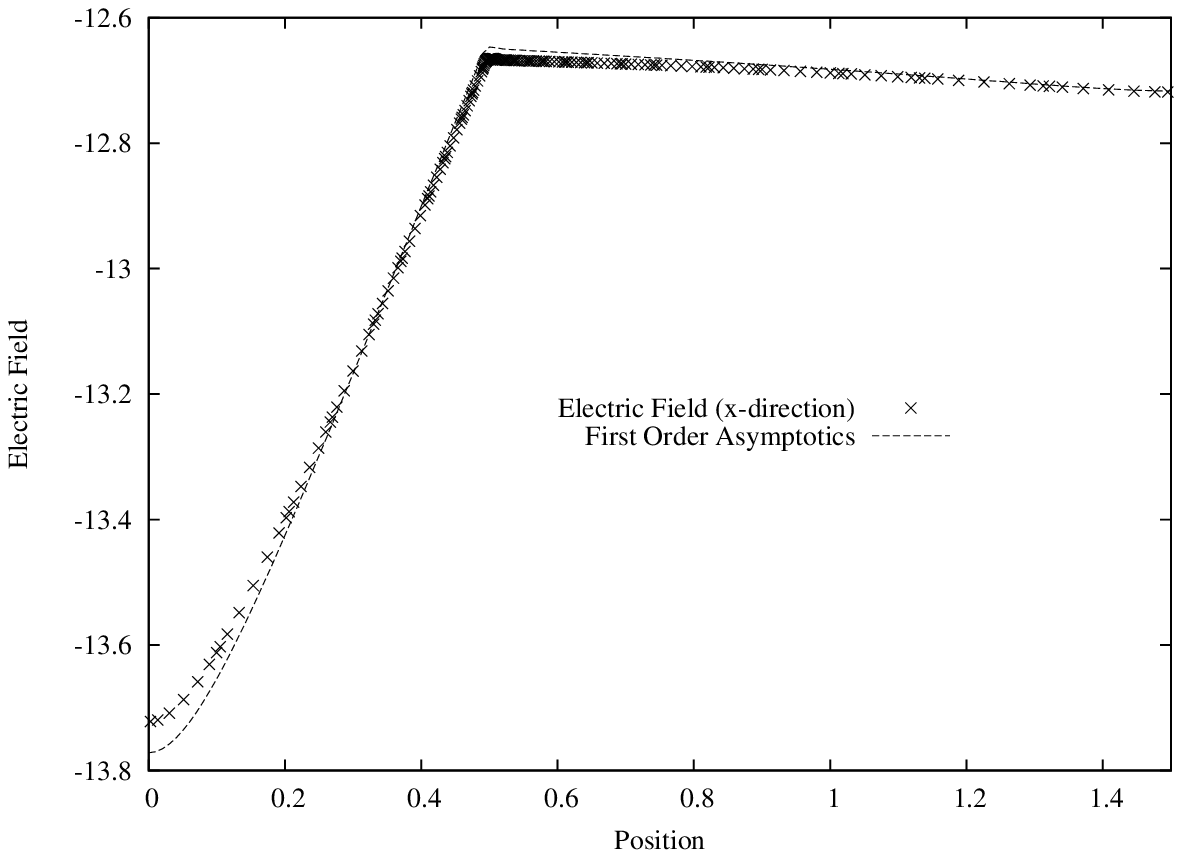}}
\caption{Plots of $n$, $p$, $X$ and $E$ for the short circuit case $V_{\mathrm{diff}} = -19.3$. Note the different scales for $n$ and $p$, and that the total variation of $E$ is less than 10\% of its value. $J= -334.2$. The electrons and holes are each dominant in one half of the device, and the asymptotic results provide good agreement to the numerical results in all cases.} \label{fig:SC}
\end{figure}

Since short circuit considers the largest negative potential difference it is thus the case for which the zero-order asymptotics derived in Sec. \ref{sec:asymptoticsinV} are best suited. This is confirmed by Fig. \ref{fig:SC} which compares the numerical solutions of $n$, $p$, $X$ and $E$ with the corresponding asymptotics. Note that we see a dramatic change in $n$, $p$, $X$ at the interface, and that either electrons or holes are largely dominant on one side of the device. Moreover, note that the slope of the electric field is significantly larger on the left-side of the device due to the larger concentration of holes on the left-side compared to the smaller concentration of electrons on the right-side (Fig. \ref{fig:SC} plots $n$ and $p$ on different scales).

The numerically calculated current (as listed in Fig. \ref{fig:SC}) is $J=-334.2$. The current predicted by \eqref{eq:JnParameters}, \eqref{eq:JpParameters} is $J^0 = -345.15$. Although not in exact agreement with the numerically calculated current for the short circuit case, this constitutes a reasonable match for the zeroth order term. Furthermore, the lowest order current (in $\delta$, $\eta$) predicted by \eqref{eq:JApprox} is $J_{\mathrm{approx}} = -339.32$, also in excellent agreement for the short circuit case. However, the current given by the asymptotics does not depend on the potential as dramatically as the simulations suggest, most likely indicating that the recombination term becomes increasingly important as we move away from short circuit.

\begin{figure}[htbp] 
\centering
\scalebox{0.8}{\includegraphics{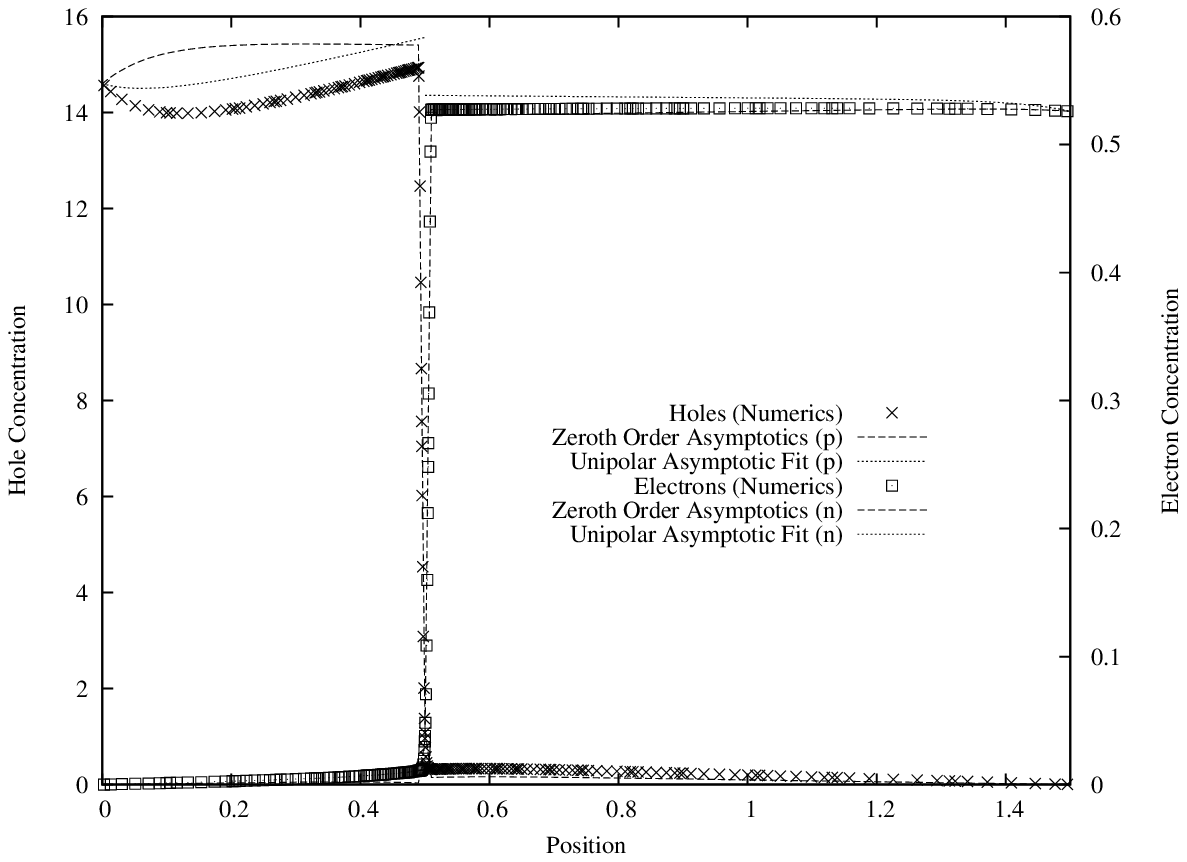}}
\caption{Plot of $n$ and $p$ for the short circuit case with artificial boundary conditions, $n_L=0.526$ and $p_0=15.6$. $J=-332.528$, compare these plots to those for the correct boundary conditions in Fig. \ref{fig:SC}. The deviation of $p$ from the constant solution results from the self-consistent potential. The effect is less pronounced for the electrons because the concentration is much lower.} \label{fig:SCBL}
\end{figure}

Another observed feature of Fig. \ref{fig:SC} are the boundary layers for $n$ to match the boundary data at $x_L$ and for $p$ to match the boundary data at $x_0$. The values of $n$ and $p$ at their majority side of the interface are predicted by the asymptotics to be $n(x_r)=0.526$ and $p(x_l)=15.6$, in good agreement with the numerically calculated values. We emphasize that Eq. \eqref{eq:JApprox} in Sec. \ref{sec:asymptoticsinV} predicts that the concentrations of $n$ and $p$ at the interface are predominantly a consequence of the leading order flux-changes over the interface $J^0_{\mathrm{approx}} = - \int_{x_l}^{x_r}{k_{d,in} X^0(y) dy}$  and are essentially independent of the boundary conditions $p_0$ and $n_L$. We can confirm this statement by using the interface values as the boundary conditions (i.e. $n(x_L) = 0.526$ and $p(x_0)=15.6$). We plot the results in Fig. \ref{fig:SCBL} without visible difference in the concentrations of the charge carriers at the interface compared to Fig. \ref{fig:SC}. In particular note that the current and the maximal values of $n$ and $p$ do not change significantly.

\subsubsection*{Optimal Power Point}
\begin{figure}[htbp] 
\centering
\scalebox{0.8}{\includegraphics{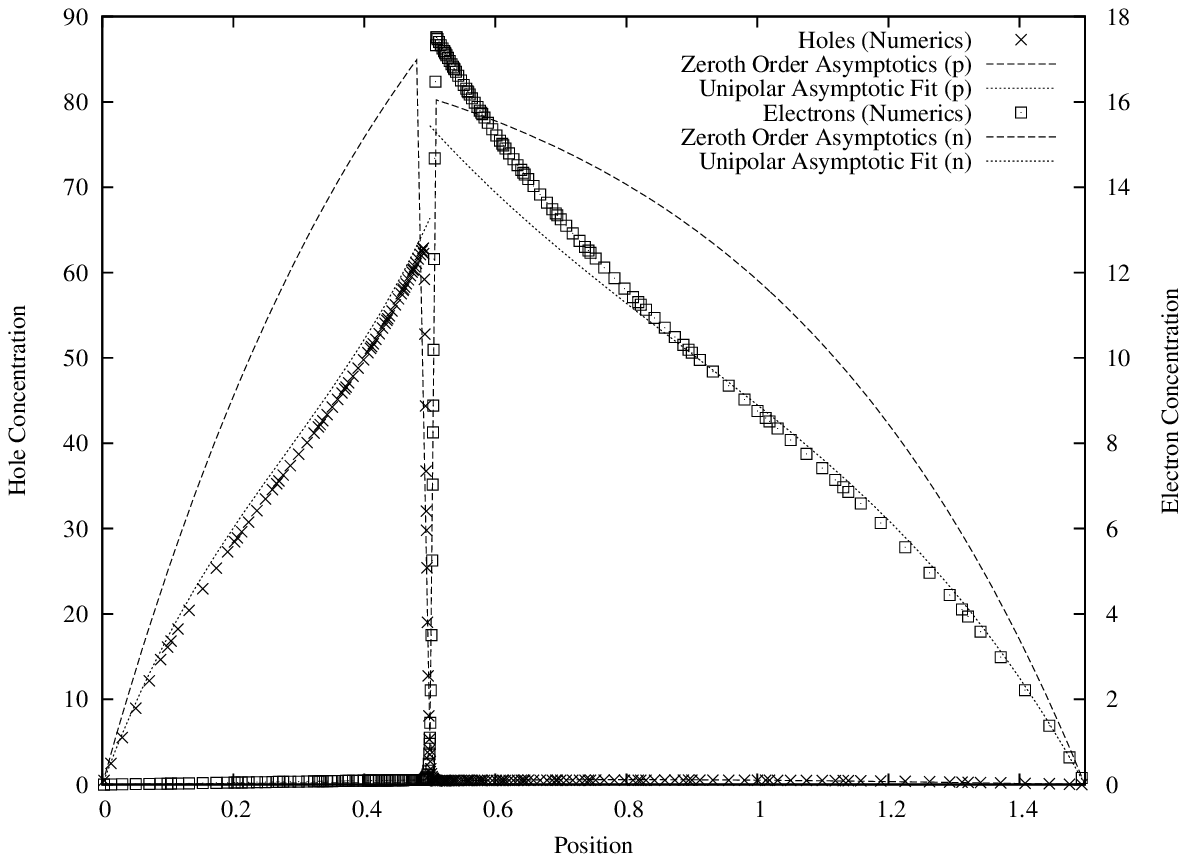}}
\scalebox{0.5}{ \includegraphics{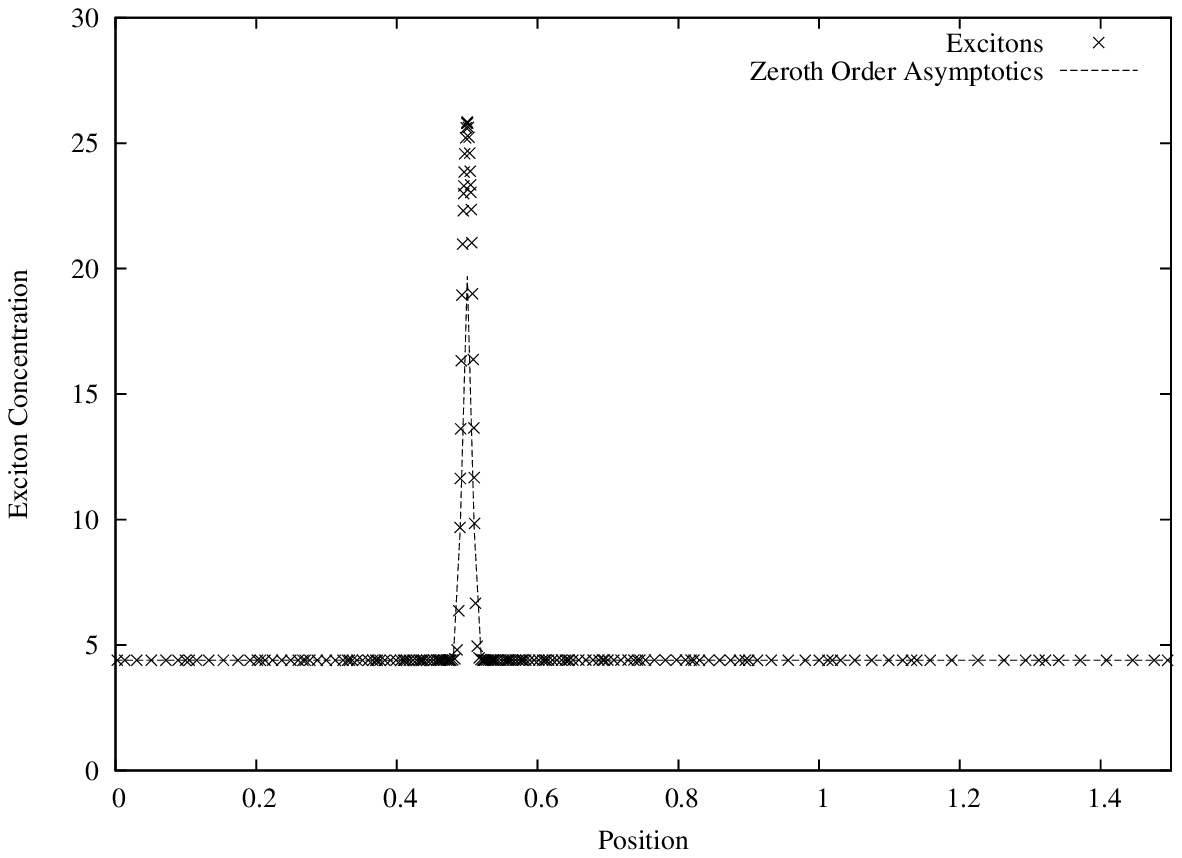} \includegraphics{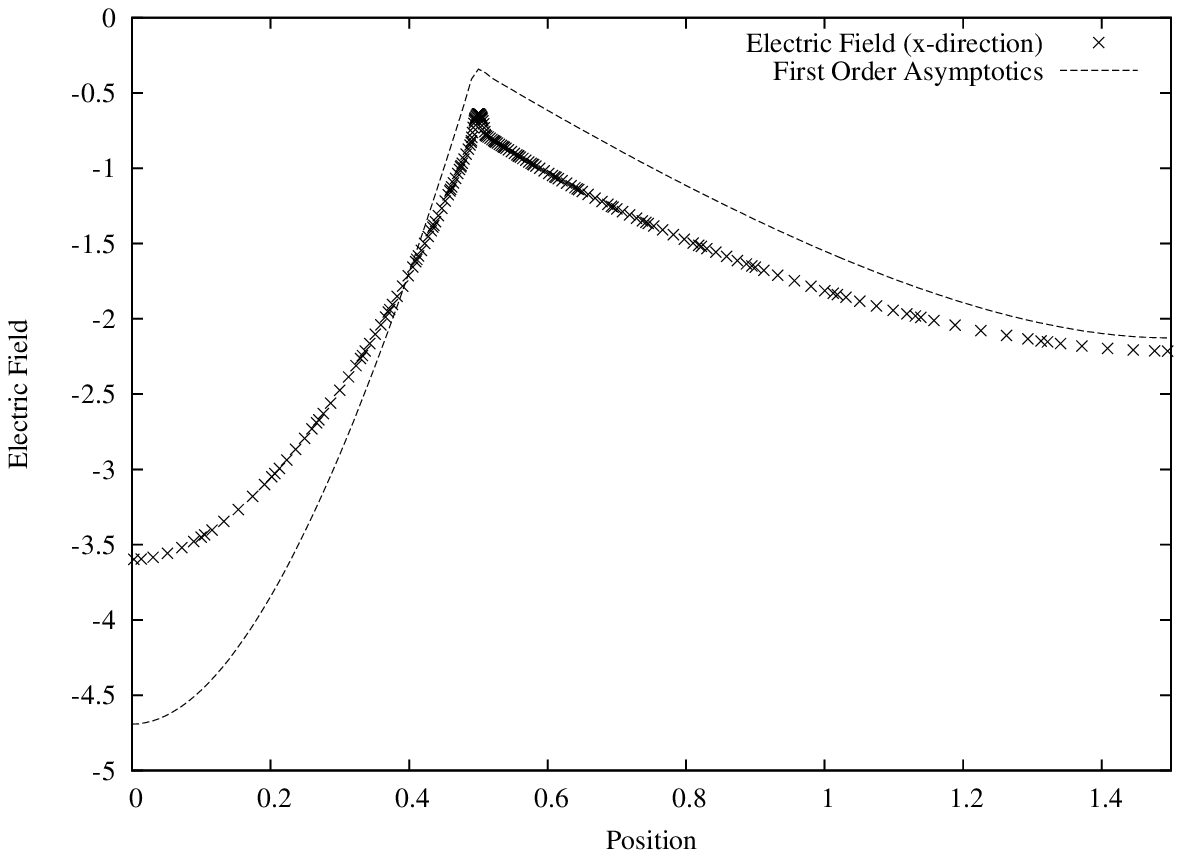}}
\caption{Plots of $n$, $p$, and $E$ respectively for the optimal power point OPP $V_{\mathrm{diff}}=-3$. The resulting current is $J = -241.585$. The exciton contribution from Poisson's equation is visible as a bump in the field plot. We see that the zeroth order asymptotics begin to diverge from the numerical and unipolar solutions as the potential difference increases.} \label{fig:OPP}
\end{figure}

A second key characteristic point for device operation is the optimal power point. This is the point where the power from the device is maximized. See Fig. \ref{fig:IVBCs} for a plot of this quantity. As shown above, this occurs for a voltage difference of $V_{\mathrm{diff}} \approx -3$. This is just below the voltage for which the resulting self-consistent field touches zero from below and the field-dependence of $k_d(E)$ causes a kink in the IV curve. The concentrations $n$, $p$ and $X$ are plotted along with the $E$-field in Fig. \ref{fig:OPP}. In the plot of the $E$-field we see a slightly more  pronounced effect (compared to SC) of the excitons on the electric field (the small bump located at the interface), but again this effect is small compared with the overall changes in the field. 

Comparing Fig. \ref{fig:OPP} to Fig. \ref{fig:SC} we see that $n$ and $p$ again are primarily located on their specific sides of the device ($p$ to the left, $n$ to the right) with an even larger change in concentrations across the interface. In contrast to Fig. \ref{fig:SC} (which features boundary layers for $n$ and $p$), the concentrations of $n$ and $p$ in Fig. \ref{fig:OPP} appear approximately linear in their majority sides. It is not clear whether this approximate linearity is a characteristic of the optimal power point or coincidental. Our numerical simulations, however, seem to suggest that this is indeed the case and we may conjecture that it is caused by a certain balance between drift and diffusion terms. Moreover, note that the scales for $n$ and $p$ remain different, but are considerably more balanced than in Fig. \ref{fig:SC}. 

Fig. \ref{fig:OPP} also plots the zeroth order asymptotics, which are still reasonable. Moreover, we see excellent fits from the unipolar asymptotic behavior of the quasi-linear behavior of the $n$ and $p$. The increasing discrepancy of the zeroth order asymptotics suggests the increasing importance of the nonlinearity of the self-consistent $E$-field. This can be confirmed by using $E^0$ instead of the numerical values $E(x_0)$ and $E(x_L)$ in the unipolar approximation which is then no longer nearly as accurate.

\subsubsection*{Open Circuit}
The third point of importance for device operation is the open-circuit voltage. This is the point at which the current through the device is zero (analogous to an open-circuit in electronics). For our device this occurs at approximately $V_{\mathrm{diff}}\approx12$. In physical units this corresponds to a potential difference of $0.31$ Volts. Considering the internal voltage of our typical device, this corresponds to an open-circuit voltage of $0.81$ Volts, larger than expected for a typical OPV. 

\begin{figure}[htbp] 
\centering
\scalebox{0.8}{\includegraphics{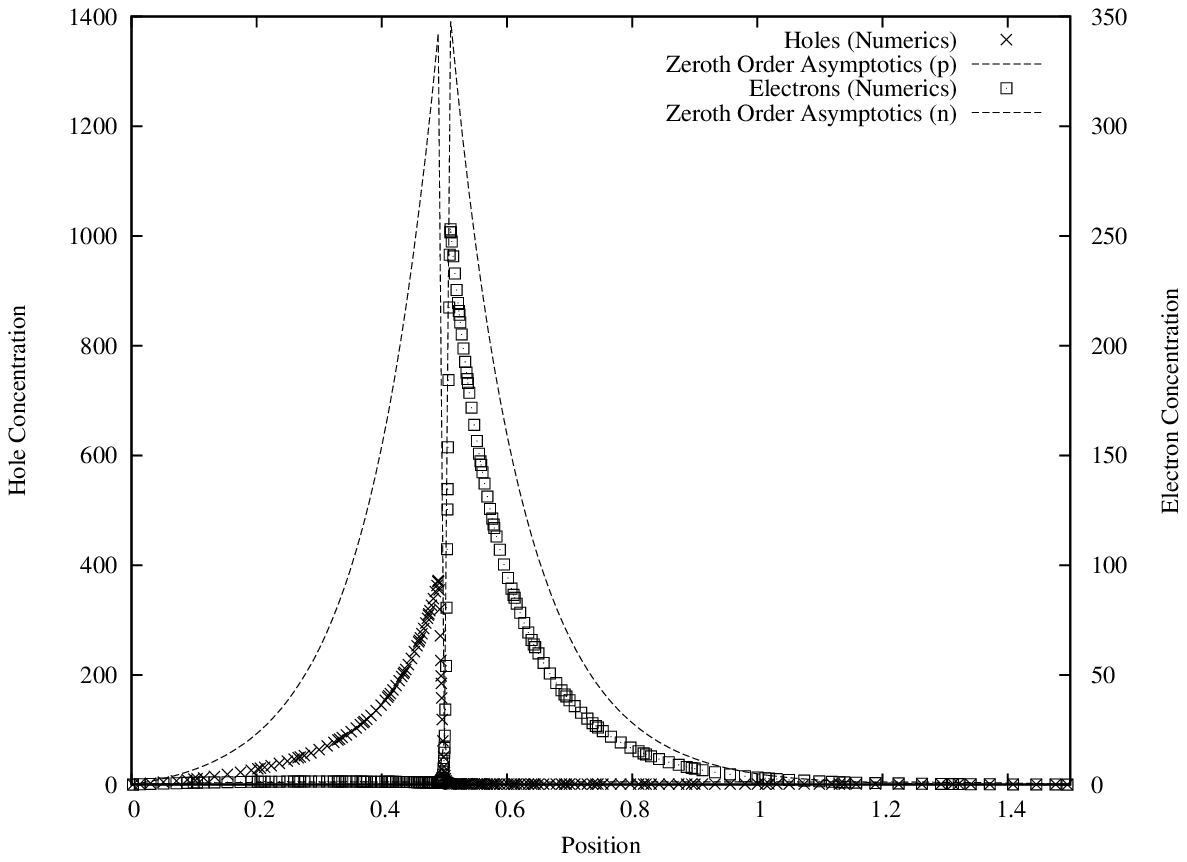}}
\scalebox{0.5}{ \includegraphics{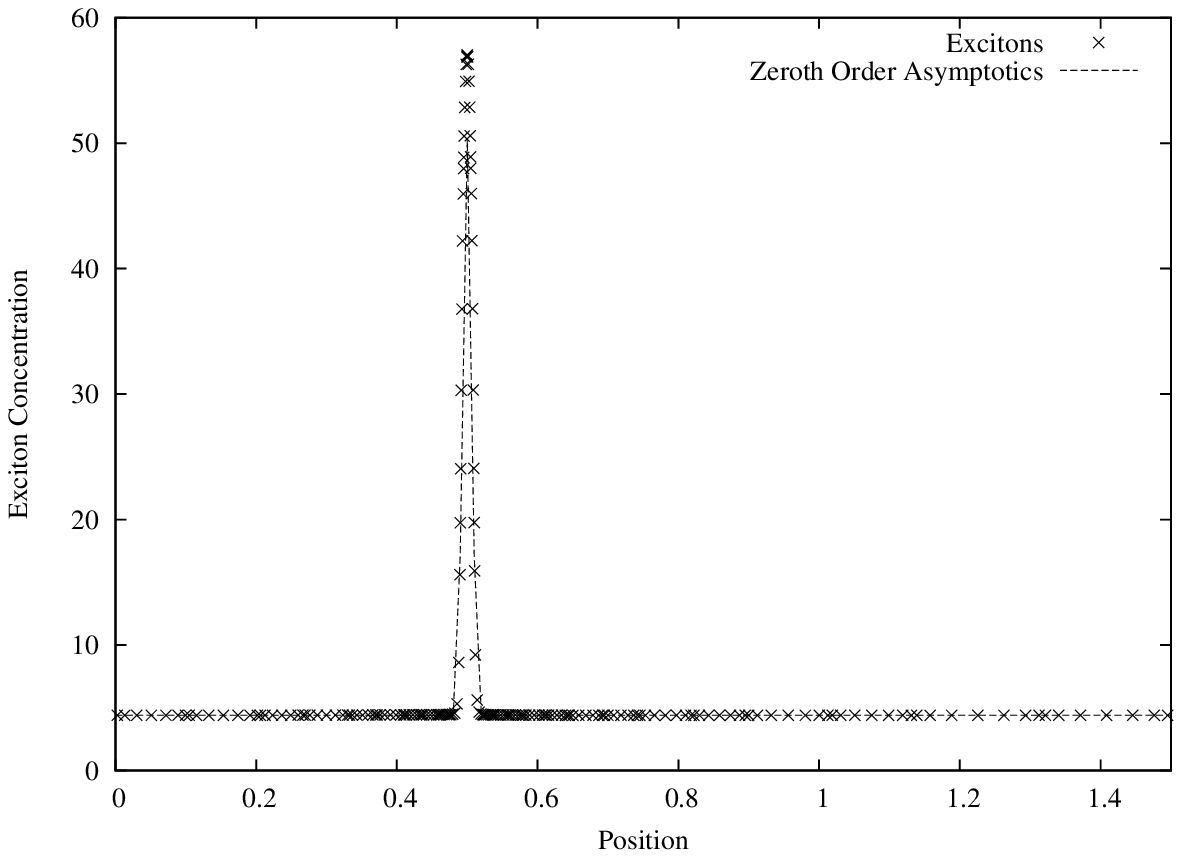} \includegraphics{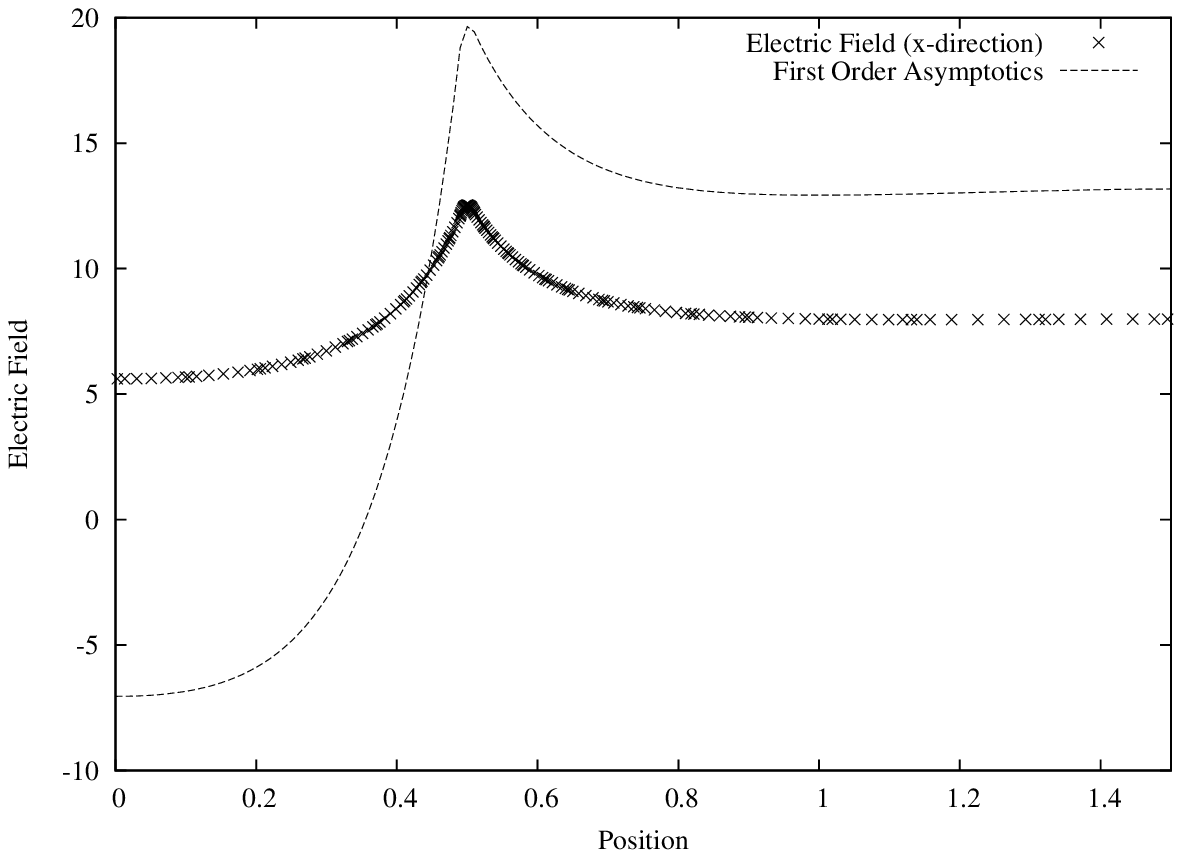}}
\caption{Plots of $n$, $p$, and $E$ respectively for the (nearly) Open-Circuit $V_{\mathrm{diff}}=12$ case.  $J = 0.23$. (The true Open Circuit Voltage is within $0.05$ of $12$, but further precision is not numerically relevant). Note the different scales for $n$ and $p$. The qualitative behavior of the zeroth order asymptotics is nearly correct, but they no longer closely match the numerical results. In particular, the estimated value for $p$ is off by a factor of $2$ nearly everywhere as the majority carrier.} \label{fig:OC}
\end{figure}

The plots in Fig. \ref{fig:OC} show a strong alteration of the concentrations $n$ and $p$ compared to the previous Fig. \ref{fig:SC} and \ref{fig:OPP}.  The hole concentration, $p$, changes from $p(x_l) \approx 370$ to $p(x_r) \approx 1.4$, over the interface, confirming that our numerical scheme is capable of capturing large changes across the interface (predicted to be $O(10^3)$ by the zeroth order asymptotics). Note that for this case the maximum values of $n$ and $p$ are of about the same order albeit plotted on different scales (by a factor four). 

The behavior of the $E$-field throughout the device is dominated by the large values of $n$ and $p$ with very little effect of the excitons at the interface. Note that we do not see a very good match of the asymptotic electric field, which is to be expected. In particular, since the zeroth asymptotics ignore the quadratic recombination term ($c_r n p$), we are neglecting this loss term for $n$ and $p$. This is particularly important within the interface, where the concentrations intersect at $n \approx p \approx 18$ (recall that $p$ is strongly decreasing from the left of the interface to the right while $n$ is strongly increasing). This leads $n$ and $p$ to be overestimated. Note, however, that the asymptotics correctly reproduce the qualitative behavior of the carrier concentrations with exponential growth from the boundaries of the device to the interface.

Finally we remark that the presented numerical scheme is capable of continuing further into forward-bias. However, the above model and in particular the parameters are primarily focused on the working-case, and we do not necessarily expect good results for the strong forward-bias regime. 
In particular, as the number of polaron pairs on the interface increases (due both to a reduction in $k_d$ and an increase in $c_r' n p$), we observe that some of them diffuse away from the interface. This corresponds physically to polaron pairs becoming excitons, a transformation which is not generally physically observed.

\subsection{Shunt-resistance and Asymptotics of the IV Curves}
\begin{figure}[htbp] 
\centering
\scalebox{.8}{\includegraphics{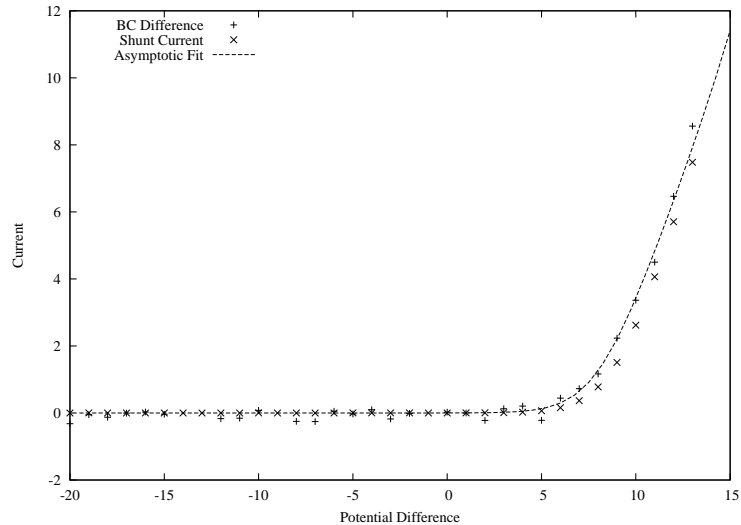}}
\caption{An IV-curve showing the similarity of the shunt current produced by injected carriers with no light generation ($\times$) and the result of subtracting the IV curve for zero boundary conditions from the usual IV curve ($+$). See Fig. \ref{fig:IVBCs} for individual plots of these curves. The dotted line represents the zeroth order asymptotic fit for the current with $G=0$. Note that the apparent noise for negative fields is a result of rounding errors in taking the difference of the currents from two different simulations.} \label{fig:IVLeak}
\end{figure}

A comparison of the IV curves in Fig. \ref{fig:IVBCs} indicates clearly that the nonlinear boundary conditions proposed in Ref. \refcite{SM} appear to have very little effect on the IV-curve. The resulting difference is barely noticeable until we reach open-circuit after which it does grow as we pass into the forward bias regime. We can study this discrepancy more clearly by simulating the OPV device in the dark with $G=0$. The resulting current is plotted in Fig. \ref{fig:IVLeak} and represents the shunt current in the device. Recalling the asymptotic currents Eq. \eqref{eq:JnParameters} and \eqref{eq:JpParameters}, we can see the increase in the shunt resistance for $E \gg 0$ (and thus $\delta \gg 1$) proportional to the terms containing $n_L$ and $p_0$. In the large $\delta$ limit we recover a term of the form: $(n_L+p_0)E$. This has the usual form of a resistor following Ohm's law with (scaled) resistivity $1/(n_L+p_0)$. 

The equivalence of this term to the difference between the nonlinear Dirichlet and homogeneous Dirichlet boundary conditions is heuristically clear from the curves plotted in Fig. \ref{fig:IVLeak}. The difference shows exactly the effect of the boundary conditions given by Ref. \refcite{SM} compared to zero B.C. For a device in the working case the boundary conditions are relatively unimportant, but as we move into the forward bias regime, these boundary conditions become more important. The given shunt resistance is reminiscent of the exponential growth of the dark-current for a photovoltaic device in forward bias. Although both currents result from taking $G=0$, our model is not designed to reproduce this dark current, which is normally generated by an equilibrium concentration term in the recombination rate - for instance $c_r ( n p - n_\infty p_\infty)$. With light generation, this term is negligible compared to the exciton contribution and thus is not included in our modeling assumptions. However, carefully choosing boundary conditions would allow our model to replicate the physically observed diode behavior in the forward bias regime. 

Note that the apparent noise for negative fields is a result of rounding errors in taking the difference of the currents from two different simulations. Choosing a smaller tolerance parameter for the convergence or calculating the current to higher precision (to eliminate $O(.1)$ errors for $O(100)$ currents) will reduce the apparent noise.

\begin{figure}[htbp] 
\centering
\scalebox{0.6}{\includegraphics{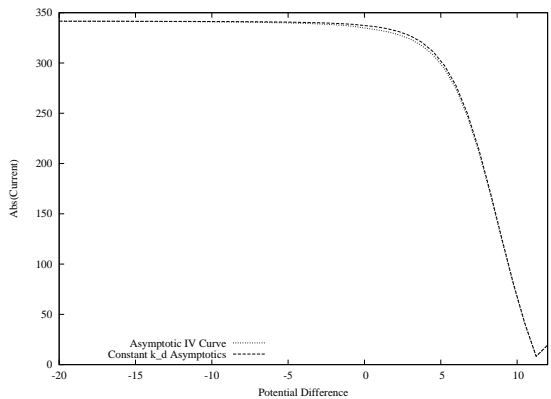}}
\caption{The IV Curve predicted by the asymptotics in Sec. \ref{sec:asymptoticsinV}. Note that it has a similar open circuit voltage and short circuit current, but has a better fill factor than the numerical simulations.} \label{fig:AsymptoticIV}
\end{figure}

Next, Fig. \ref{fig:AsymptoticIV} plots the asymptotic IV curve predicted by the sum $J_{n0}+J_{p0}$ of Eqs. \eqref{eq:JnParameters} and \eqref{eq:JpParameters}. Although the asymptotic approximation is questionable for small potential differences (as stated previously) we do observe a qualitatively realistic IV-curve in Fig. \ref{fig:AsymptoticIV}. In particular we observe a similar open-circuit voltage and short circuit current as well as a significantly better fill-factor of 87\%.
The high fill factor for the asymptotic curves is expected, since the asymptotics neglect the recombination current, which  should be the main cause of power loss  for a bilayer device. Indeed, by ignoring the recombination term, the asymptotics describe what is called the photogenerated current, the idealized current generated by the device. This has interesting properties in and of itself, and this avenue of prediction is under further investigation.

\begin{figure}[htbp] 
\centering
\scalebox{0.6}{\includegraphics{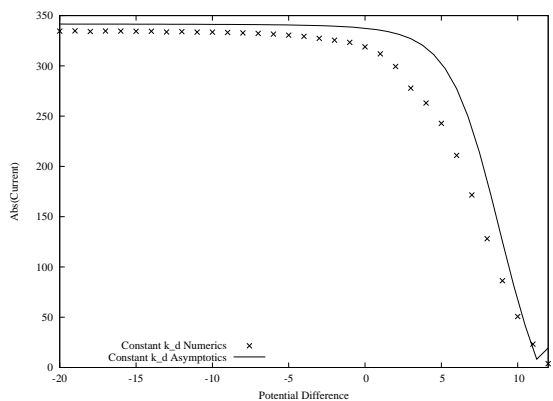}}
\caption{A comparison of the IV-curves obtained using the numerical and asymptotic schemes with $k_{d,in} =k_{d,in}(E_{SC})= 2763$ to be constant.} \label{fig:IVkdconst}
\end{figure}

A comparable increase in fill-factor is observed in Fig. \ref{fig:IVkdconst}, where we plot a numerical IV curve assuming a constant dissociation rate $k_{d,in} \gg k_{r,in}$. Moreover, Fig. \ref{fig:IVkdconst} plots the corresponding asymptotic IV curve with good agreement. This indicates that a primary loss of efficiency is due to the exciton recombination term, which becomes non-negligible near the interface when approaching the open-circuit voltage. This quantifies the physical intuition that device efficiency is increased by forcing the polaron pairs to dissociate more efficiently into free carriers. For the numerical simulation the new fill factor is a much improved 60\%, whereas the asymptotics are virtually unchanged from the standard model shown in Fig. \ref{fig:AsymptoticIV}.

We  therefore plan to further study which parameter changes give the most improvement in asymptotic device characteristics. By studying these various parameters, we shall aim to determine which devices are mostly likely to have the optimal characteristics, and then test whether they also imply an improved efficiency for the full numerical system. We expect to find interesting information in terms of device design and implementation on a fairly general level rather than optimizing our model to reproduce one particular pair of organic materials. Further work is underway to establish whether the predictions created by our model are directly applicable from a device design standpoint. The IV-plots shown here match qualitatively to those of realistic OPVs. A precise quantitative agreement will be challenging due to the large number of parameters in the problem, many of which are difficult to measure using standard physical techniques.

Other further work is in progress addressing the rigorous mathematical theory of the system presented, as well as the possibility of deriving similar models for which the size of the interface tends to zero.

\begin{figure}[htbp] 
\centering
\scalebox{0.3}{\includegraphics{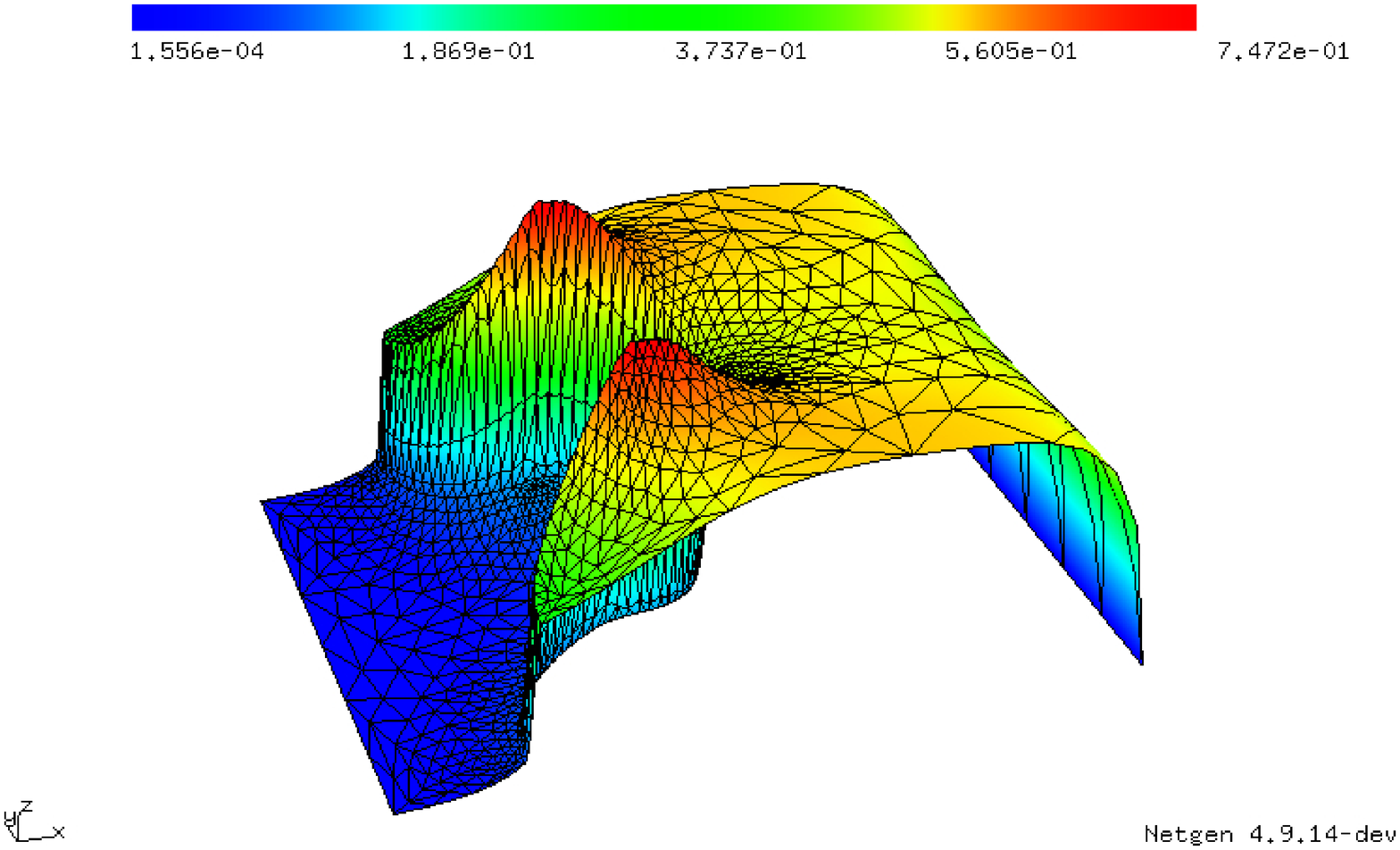}}
\scalebox{.3} {\includegraphics{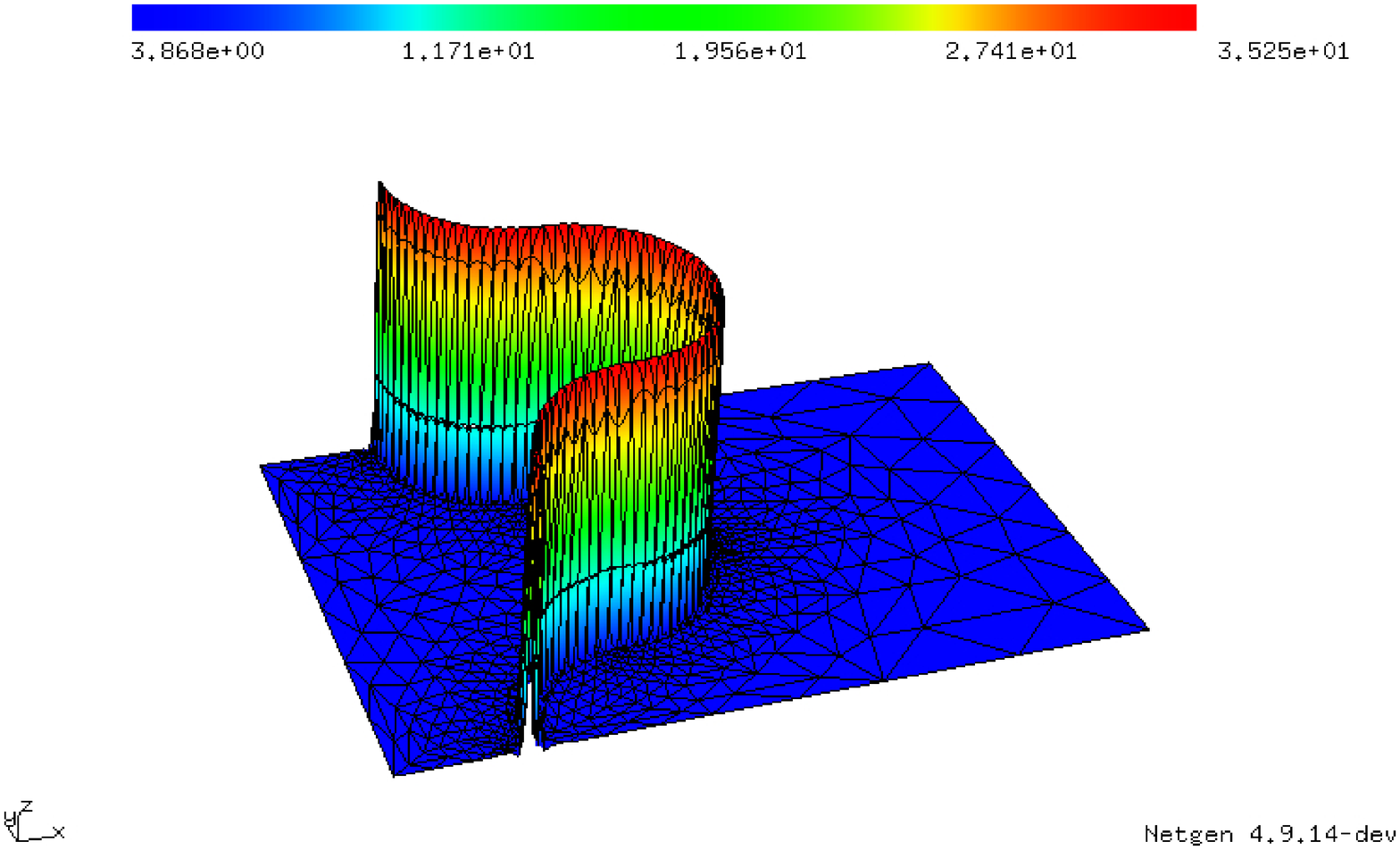}}
\caption{Sample numerical results for the electron and hole concentrations for a device with a non-uniform interface at short-circuit. The exciton concentration shows good confinement and demonstrates the shape of the interface. The strange peaks in $n$ appear as a result of the given $C^1$ form of $U$ and the concentrating effects of the electric field in the vertical direction. Additional modelling assumptions will be necessary to obtain physically realistic results. Note that away from the interface we see the same nearly-constant behavior with a boundary layer near the metal contact.} \label{fig:2Dplot}
\end{figure}

Current interest also revolves around the development of intricate interface geometries in solar cells. The most efficient solar cells do not use simple bilayer interfaces, but generally use heterojunction interfaces which maximize the active area of the cell. Our model can be extended to the case of a bilayer with regular periodic interface, with the primary difficulty being the specific modeling assumptions for non-straight interfaces. In Fig. \ref{fig:2Dplot} we show sample plots of the electrons and excitons for a non-uniform interface to demonstrate the capabilities of the numerics.

\appendix
\section{Calculation of $\Phi_n(x)$ and  $\Phi_p(x)$} \label{ap:phi}
First note that $\varphi_x(x) = -E^0$ outside the interface and $\varphi_x(x) = -E^0 + U_x(x)$ inside the interface. In addition $\varphi(b)-\varphi(a) = -E^0(b-a)$ outside the interface and $\varphi(b)-\varphi(a)=-E^0(b-a) + U(b)-U(a)$ inside the interface. Note that the change in form of $\varphi$ forces different treatment of $\varphi(y)$ depending on the region of integration and a change in $\varphi(x)$ depending on the value of $x$. We proceed by cases:
\begin{align*}
 \Phi_n&=\begin{cases}
    \int_{x_0}^x{e^{-E^0 (x-y)}dy} \qquad &x_0 < x \leq x_l \\[1mm]
    \int_{x_0}^{x_l}{e^{-E^0 (x-y) + U(x)-U(x_l)}dy} +  \int_{x_l}^x{e^{-E^0(x-y) + U(x)-U(y)}dy} \qquad &x_l < x \leq x_r \\[1mm]
    \int_{x_0}^{x_l}{e^{-E^0 (x-y) + U(x_r)-U(x_l)}dy} + \\
    \qquad +\int_{x_l}^{x_r}{e^{-E^0(x-x_l) + U(x_r)-U(y)}dy} + \int_{x_r}^x{e^{-E^0(x-y)}dy} \qquad &x_r < x < x_L
\end{cases} \\[3mm]
&= \begin{cases}
    \frac{1}{E^0} \left(1 - e^{-E^0 (x-x_0)}\right)   \qquad \quad &x_0 < x \leq x_l \\[1mm]
    e^{-E^0(x-x_l) + U(x)-U(x_l)}\Phi_n(x_l) + \frac{1 - e^{-E^0(x-x_l) + U(x)-U(x_l)}}{E^0-U_x(\theta)}   \qquad \quad &x_l < x \leq x_r \\[1mm]
     e^{-E^0(x-x_r)}\Phi_n(x_r) + \frac{1}{E^0} \left(1 - e^{-E^0(x-x_r)}\right)   \qquad \quad &x_r < x < x_L
\end{cases}
\end{align*}
where we have included the $U(x_l)$ terms for clarity even though usually we take $U(x_l)=0$. In addition, the term $U_x(\theta)$ for $\theta \in [x_l, x_r]$ arises from using the mean value theorem after integrating by parts once. If we assume that $U$ is piecewise linear, then $U_x = \frac{U(x_r)-U(x_l)}{x_r-x_l}$ is constant inside the interface and the integral is straightforward. In the same manner, we have:
\begin{equation*}
 \Phi_p = \begin{cases}
    \frac{-1}{E^0} \left(1 - e^{E^0 (x-x_0)}\right) \qquad &x_0 < x \leq x_l \\[1mm]
    e^{E^0(x-x_l) - U(x)+U(x_l)}\Phi_p(x_l) - \frac{1 - e^{E^0(x-x_l) - U(x)+U(x_l)}}{E^0-U_x(\theta)}  \qquad &x_l < x \leq x_r \\[1mm]
     e^{E^0(x-x_r)}\Phi_p(x_r) - \frac{1}{E^0} \left(1 - e^{E^0(x-x_r)}\right) \qquad &x_r < x < x_L
\end{cases} \text{.}
\end{equation*}

Recalling the definitions given in Eq. \ref{eq:deltaeta}, $\eta=e^{E^0(x_r-x_l) - U(x_r)+U(x_l)}$, $\delta= e^{E^0(x_l-x_0)}$ as well as $\delta^2=e^{E^0(x_L-x_r)}$, we can write the particular values $\Phi_n(x_L)$ and $\Phi_p(x_L)$ simply as
\begin{align*}
 \Phi_n(x_L)&=\frac{1}{E^0} \frac{1}{\eta \delta^2}\left(1 - \frac{1}{\delta}\right) + \frac{1}{E^0-U_x(\theta)} \frac{1}{\delta^2} \left(1 - \frac{1}{\eta}\right) + \frac{1}{E^0} \left(1 - \frac{1}{\delta^2}\right) \\
\Phi_p(x_L)&=\frac{-1}{E^0} \eta \delta^2 \left(1 - \delta \right) - \frac{1}{E^0-U_x(\theta)} \delta^2\left(1 - \eta\right) - \frac{1}{E^0} \left(1 -\delta^2\right) \text{.}
\end{align*}
Note that the for $E^0<0$, we have $\delta \ll 1$, and the results are positive (as expected for integrating a positive quantity).

\section*{Acknowledgment}
The authors acknowledge support from King Abdullah University of Science and Technology (KAUST) Award Number: KUK-I1-007-43. PAM also acknowledges support from the Fondation Sciences Mathematique de Paris, in form of his Excellence Chair 2011/2012. MTW acknowledges financial support from the Austrain Science Foundation (FWF) via the Hertha Firnberg project TU56-N23.
KF acknowledges the support of NaWi Graz. 

\bibliography{PhotovoltaicReferences}{}
\bibliographystyle{plain}
\end{document}